\documentclass[12pt,a4paper]{amsart}

\usepackage{graphicx}

\usepackage{amssymb, subfig, amsmath, mathrsfs, mathtools}
\usepackage[a4paper, margin = 1in]{geometry}
\usepackage{microtype, fancyhdr, polynom, amsthm}
\usepackage{amscd, color, tikz, tikz-cd, array}

\usepackage[hyperindex,breaklinks]{hyperref}

\usepackage{amsfonts}
\usepackage[latin2]{inputenc} 
\usepackage{t1enc}
\usepackage[mathscr]{eucal}
\usepackage{indentfirst}
\usepackage{graphics}
\usepackage{pict2e}
\usepackage{epic}
\usepackage{epstopdf}

\numberwithin{equation}{section}

\usepackage[numbers,sort&compress]{natbib}

\usetikzlibrary{decorations.markings}

\newcommand{\m}{\mathbb}
\newcommand{\C}{\mathscr}
\newcommand{\B}{\mathcal}

\newcommand{\V}{\mathsf{Vec}}

\newcommand{\Aut}{\operatorname{Aut}}

\newcommand{\im}{\operatorname{Im}}

\newcommand{\Rep}{\operatorname{Rep}}

\newcommand{\iso}{\xrightarrow{\sim}}

\newcommand{\tr}{\operatorname{tr}}
\newcommand{\ov}{\overline}
\newcommand{\GL}{\operatorname{GL}}

\newcommand{\PSL}{\operatorname{PSL}}
\newcommand{\PL}{\operatorname{PGL}}

\newcommand{\Span}{\operatorname{Span}}
\newcommand{\rank}{\operatorname{rank}}

\newtheorem{mythm}{Theorem}[section]

\newtheorem{mydef}[mythm]{Definition}
\newtheorem{mylem}[mythm]{Lemma}
\newtheorem{mycor}[mythm]{Corollary}

\allowdisplaybreaks
\makeatletter
\tikzset{nomorepostaction/.code=\let\tikz@postactions\pgfutil@empty}
\makeatother

\newcommand{\KP}[1]{%
  \begin{tikzpicture}[baseline=-\dimexpr\fontdimen22\textfont2\relax]
  #1
  \end{tikzpicture}%
}

\newcommand{\BoxedHom}{}
\def\BoxedHom[#1, #2](#3, #4)(#5){

    \draw[color=gray,thick] (#1 - #3/2, #2 - #4/2) -- (#1 + #3/2, #2 - #4/2) -- (#1 + #3/2, #2 + #4/2) -- (#1 - #3/2, #2 + #4/2) -- cycle;

    \draw (#1, #2) node  {#5};
}

\newcommand{\ModularDataSDetailed}{%
    \KP{%
        \begin{scope}[decoration={
        markings,
        mark=at position 0.5 with {\arrow{>}}}
        ]
          \draw[color=gray, thick, postaction={decorate}] (0, 0) -- (0, 7.6);
          \draw[color=gray, thick, postaction={decorate}] (6, 9) -- (6, 0.5);
        \end{scope}

        \draw[color=gray, thick] (3, 0) arc (0:-180:1.5);
        \draw[color=gray, thick] (3, 0) -- (3, 2.5) -- (4, 3.5) -- (4, 4.5) -- (3.7, 4.8);

        \draw[color=gray, thick] (6, 0.5) arc (0:-180:0.5);

        \draw[color=gray, thick] (5, 0.5) -- (5, 1.5) -- (4, 2) -- (4, 2.5) -- (3.7, 2.8);

        \draw[color=gray, thick] (3.3, 3.2) -- (3, 3.5) -- (3, 4.5) -- (4, 5.5) -- (4, 6) -- (5, 6.5) -- (5, 7.6);

        \draw[color=gray, thick] (3.3, 5.2) -- (3, 5.5) -- (3, 9.5);

        \BoxedHom[5, 8](0.8, 0.8)($\phi_W$);

        \BoxedHom[0, 8](0.8, 0.8)($\phi_V$);

        \draw[color=gray, thick] (0, 8.4) -- (0, 9.5);

        \draw[color=gray, thick] (5, 8.4) -- (5, 9);

        \draw[color=gray, thick] (6, 9) arc (0:180:0.5);

        \draw[color=gray, thick] (3, 9.5) arc (0:180:1.5);

        \draw (12, -2) node {$I$};
        \draw (12, -0.75) node {$V \otimes (V^* \otimes I)$};
        \draw (12, 0.875) node {$V \otimes \big(V^* \otimes (W \otimes W^*)\big)$};
        \draw (12, 2.375) node {$V \otimes \big((V^* \otimes W) \otimes W^*\big)$};
        \draw (12, 4) node {$V \otimes \big((W \otimes V^*) \otimes W^*\big)$};
        \draw (12, 5.675) node {$V \otimes \big((V^* \otimes W) \otimes W^*\big)$};
        \draw (12, 7.125) node {$V \otimes \big(V^* \otimes (W \otimes W^*)\big)$};
        \draw (12, 8.75) node {$V^{**} \otimes \big(V^* \otimes (W^{**} \otimes W^*)\big)$};
        \draw (12, 10.25) node {$V^{**} \otimes (V^* \otimes I)$};
        \draw (12, 11.5) node {$I$};

        \draw[color=gray, dotted] (-0.5, -1.5) -- (10, -1.5);
        \draw[color=gray, dotted] (-0.5, 0) -- (10, 0);
        \draw[color=gray, dotted] (-0.5, 1.75) -- (10, 1.75);
        \draw[color=gray, dotted] (-0.5, 3) -- (10, 3);
        \draw[color=gray, dotted] (-0.5, 5) -- (10, 5);
        \draw[color=gray, dotted] (-0.5, 6.25) -- (10, 6.25);
        \draw[color=gray, dotted] (-0.5, 8) -- (-0.4, 8);
        \draw[color=gray, dotted] (0.4, 8) -- (4.6, 8);
        \draw[color=gray, dotted] (5.4, 8) -- (10, 8);
        \draw[color=gray, dotted] (-0.5, 9.5) -- (10, 9.5);
        \draw[color=gray, dotted] (-0.5, 11) -- (10, 11);

        \draw (7.6, -1.25) node {$\eta_V$};
        \draw (7.6, 0.25) node {$1 \otimes 1 \otimes \eta_W$};
        \draw (7.6, 2) node {$1 \otimes (\omega^{-1}_{V^*, W, W^*})^{-1}$};
        \draw (7.6, 3.25) node {$1 \otimes \sigma_{V^*, W} \otimes 1$};
        \draw (7.6, 5.25) node {$1 \otimes \sigma_{W, V^*} \otimes 1$};
        \draw (7.6, 6.5) node {$1 \otimes \omega^{-1}_{V^*, W, W^*}$};
        \draw (7.6, 8.25) node {$\phi_V \otimes 1 \otimes \phi_W \otimes 1$};
        \draw (7.6, 9.75) node {$1 \otimes 1 \otimes \epsilon_{W^*}$};;
        \draw (7.6, 11.25) node {$\epsilon_{V^*}$};
    }
}
\newcommand{\ModularDataTDetailed}{%
    \KP{%
        \begin{scope}[decoration={
        markings,
        mark=at position 0.5 with {\arrow{>}}}
        ]
          \draw[color=gray, thick, postaction={decorate}] (0, 6) -- (0, 0);
          \draw[color=gray, thick, postaction={decorate}] (6, 5.5) -- (6, 0.5);
        \end{scope}

        \draw[color=gray, thick] (3, 0) arc (0:-180:1.5);
        \draw[color=gray, thick] (3, 0) -- (3, 2.5) -- (4, 3.5) -- (4, 4) -- (5, 4.5) -- (5, 5.5);

        \draw[color=gray, thick] (6, 0.5) arc (0:-180:0.5);

        \draw[color=gray, thick] (5, 0.5) -- (5, 1.5) -- (4, 2) -- (4, 2.5) -- (3.7, 2.8);

        \draw[color=gray, thick] (3.3, 3.2) -- (3, 3.5) -- (3, 6);

        \draw[color=gray, thick] (6, 5.5) arc (0:180:0.5);

        \draw[color=gray, thick] (3, 6) arc (0:180:1.5);

        \draw (12, -2) node {$I$};
        \draw (12, -0.75) node {$V^* \otimes (V^{**} \otimes I)$};
        \draw (12, 0.875) node {$V^* \otimes \big(V^{**} \otimes (V \otimes V^*)\big)$};
        \draw (12, 2.375) node {$V^* \otimes \big((V^{**} \otimes V) \otimes V^*\big)$};
        \draw (12, 3.625) node {$V^* \otimes \big((V \otimes V^{**}) \otimes V^*\big)$};
        \draw (12, 5.125) node {$V^* \otimes \big(V \otimes (V^{**} \otimes V^*)\big)$};
        \draw (12, 6.75) node {$V^* \otimes (V \otimes I)$};
        \draw (12, 8) node {$I$};

        \draw[color=gray, dotted] (-0.5, -1.5) -- (10, -1.5);
        \draw[color=gray, dotted] (-0.5, 0) -- (10, 0);
        \draw[color=gray, dotted] (-0.5, 1.75) -- (10, 1.75);
        \draw[color=gray, dotted] (-0.5, 3) -- (10, 3);
        \draw[color=gray, dotted] (-0.5, 4.25) -- (10, 4.25);
        \draw[color=gray, dotted] (-0.5, 6) -- (10, 6);
        \draw[color=gray, dotted] (-0.5, 7.5) -- (10, 7.5);

        \draw (7.6, -1.25) node {$\eta_{V^*}$};
        \draw (7.6, 0.25) node {$1 \otimes 1 \otimes \eta_V$};
        \draw (7.6, 2) node {$1 \otimes (\omega^{-1}_{V^{**}, V, V^*})^{-1}$};
        \draw (7.6, 3.25) node {$1 \otimes \sigma_{V^{**}, V} \otimes 1$};
        \draw (7.6, 4.5) node {$1 \otimes \omega^{-1}_{V, V^{**}, V^*}$};
        \draw (7.6, 6.25) node {$1 \otimes 1 \otimes \epsilon_{V^*}$};
        \draw (7.6, 7.75) node {$\epsilon_V$};
    }
}

\title{Computing Modular Data for Pointed Fusion Categories}

\author[Angus Gruen, Scott Morrison]{Angus Gruen, Scott Morrison}

\address{Australian National University \\ Mathematical Science Institute \\
Acton ACT 2600} 

\email{angusgruen@gmail.com}
\email{scott.morrison@anu.edu.au}

\makeindex

\begin{document}

\begin{abstract}
  \label{cha:abstract}

A formula for the modular data of $\B Z(\V^{\omega}G)$ was given by Coste, Gannon and Ruelle in \cite{coste:2000}, but without an explicit proof for arbitrary 3-cocycles. This paper presents a derivation using the representation category of the quasi Hopf algebra $D^{\omega}G$. Further, we have written code to compute this modular data for many pairs of small finite groups and $3$-cocycles. This code is optimised by taking advantage Galois symmetries of the $S$ and $T$ matrices. We have posted a database of modular data for the Drinfeld center of every Morita equivalence class of pointed fusion categories of dimension less than $64$. 

\end{abstract}

\maketitle

\section{Introduction}

Fusion categories are an important area of study due both to their frequent occurrence in category theory as well as their applications in physics. An important current area of study in condensed matter physics is topological phases of matter. One approach to studying topological phases of matter is with topological quantum field theories (TQFTs). The cobordism hypothesis \cite{lurie:2009} states that TQFTs are classified by higher categorical data and in particular the work of Douglas--Schommer-Pries--Snyder \cite{douglas:2013} identifies fully extended TQFTs in 2+1 dimensions with fusion categories.

A pointed fusion category is one in which every object is invertible. It is well known that that every pointed fusion category is equivalent to $\V^{\omega}G$ for some group $G$ and three-cocycle $\omega$. The purpose of this paper is to study the Drinfeld centers $\B{Z}(\V^{\omega} G)$ which are modular tensor categories. These modular tensor categories also arise as representation categories of the ribbon quasi-hopf algebra $D^{\omega} G$ known as the twisted Drinfeld double of a finite group. This algebra $D^{\omega} G$, which first appeared in \cite{Pasquier:1990} has been extensively studied in mathematical physics in part due to its appearance in Dijkgraaf-Witten theory \cite{Witten:1990}. In particular, in $2 + 1$ dimensions, Dijkgraaf-Witten theory corresponds under the cobordism hypothesis to $\V^{\omega} G$ and the invariant of $S^1$ is $\B{Z}(\V^{\omega} G)$.

While the categories $\V^{\omega}G$ are relatively straighforward, the centers $\B{Z}\big(\V^{\omega}G\big)$ can be difficult to pin down. Recently in \cite{Schauenburg:2017}, Mignard and Schauenberg identified infinitely many finite groups $G$, for which there are (at least) two $3$-cocycles $\omega, \omega'$ so $\B{Z}(\V^{\omega} G) \ncong \B{Z}(\V^{\omega'} G)$, but nevertheless they have equivalent modular data. The smallest examples appear at $G = \m{Z}/5\m{Z} \rtimes \m{Z}/11\m{Z}$ and $G = \m{Z}/5\m{Z} \rtimes \m{Z}/31\m{Z}$ and the corresponding sets of modular data in both of these cases can be found in our database. It is currently unknown if the $S$ and $T$ matrices are complete invariants when $|G| < 55$ and our database can point to the set of cases left to check. 

\subsection{Outline}

	Section \ref{cha:ModData} introduces the twisted quantum double of a finite group $D^{\omega}G$, and gives the equivalence between $\Rep\big(D^{\omega}G\big)$ and $\B Z\big(\V^{\omega^{-1}}G\big)$. A detailed proof of this equivalence can be found in Appendix \ref{app:app2}. This is well known to the experts but we were unable to find a source that spells this out completely. Using this equivalence, we provide a detailed derivation for the modular data of $\B Z\big(\V^{\omega^{-1}}G\big)$. The equations of the modular data were given without proof as Equations 5.23 and 5.24 in \cite{coste:2000}. Simplified cases have been proven in \cite{coquereaux:2012, coste:2000, dijkgraaf:1990, bantay:1996} but we have not found a careful derivation of the general formula and so include one here.

	Sections \ref{cha:CompData} and \ref{cha:ana} deal with the computational aspects. (Any reader who already feels sufficiently confident in the formula from \cite{coste:2000} may like to skip directly to Section \ref{cha:CompData}) Section \ref{cha:CompData} discusses our code that computes the modular data of $\B Z\big(\V^{\omega^{-1}}G\big)$ given $G$ and $\omega$. In particular it discusses the methods used to identify the simple objects of $\B Z\big(\V^{\omega^{-1}}G\big)$ and then construct the $S$ and $T$ matrices. We detail a significant optimisation which takes advantage of the Galois symmetries of the modular data. Using this code, we assemble a database of this modular data for the Drinfeld doubles of all Morita equivalence classes of pointed fusion categories with dimension at most $63$. Section \ref{cha:ana} gives some preliminary results that can be obtained from this database. In particular we confirm the result of Mignard and Schauenburg in \cite{mignard:2017} concerning the number of Morita equivalence classes of pointed fusion categories at each dimension less than $32$, improve upon their lower bound when the dimension is equal to $32$ and publish new lower bounds at each dimension between $33$ and $63$. We note that the improvement at dimension $32$ also shows that the modular data of a category is a stronger invariant than merely the Frobenius-Schur indicators and $T$-matrix.

	The database of modular data constructed is available for public use at \url{https://tqft.net/web/research/students/AngusGruen}. This database already been put to use by AnRan Chen \cite{An-Ran:2017} computing invariants for torus knots and to study examples of grafting, as in \cite{Evans:2011}. 
\section{Modular Data for the Drinfeld Centre of Twisted \texorpdfstring{$G$}{G}-graded Vector Spaces}
\label{cha:ModData}

For simplicity, we work over the field $\m C$.
The goal of this section is to give a derivation for the formula of the $S$ and $T$ matrices for $\B Z\big(\V^{\omega^{-1}}G\big)$. This will be done by constructing a quasitriangular quasi Hopf algebra $D^{\omega}G$ and using the ribbon equivalence $\B Z\big(\V^{\omega^{-1}}G\big) \cong \Rep(D^{\omega}G)$. This allows us to translate the morphisms defining the $S$ and $T$ matrices into algebraic expressions in $D^{\omega}G$ where we can simplify them. We assume without loss of generality that $\omega$ is unitary as every cocycle is cohomologous to a unitary cocycle and if $\omega$ and $\omega'$ are cohomologous then $\V^{\omega}G \overset{\otimes}{\cong} \V^{\omega'}G$.

The main result from this section is a proof the following theorem from \cite{coste:2000}.

For a group $G$ and a unitary $3$-cocycle $\omega$ define
\begin{equation}
    \theta_{g}(x, y) = \frac{\omega(g, x, y)\omega(x, y, (xy)^{-1}gxy)}{\omega(x, x^{-1}gx, y)} \label{Derived two cocycle}.
\end{equation}
Next, for each conjugacy class $K$ of $G$, fix a representative $a$ and for every $g \in K$ a group element $a_g$ such that $a = a_gga_g^{-1}$. Denote by $C(a) = \{g \in G | ag = ga\}$ the centraliser of $a$.

\begin{mythm} \label{S and T matrices}
    The simple objects of $\B Z\big(\V^{\omega^{-1}}G\big)$ correspond to tuples $(K, a, \rho_a)$ with $K$ a conjugacy class of $G$ with representative $a$ and $\rho_a$ an irreducible $\theta_a$-projective representation of $C(a)$. Then the entries of the $S$ and $T$ matrices of $\B Z\big(\V^{\omega^{-1}}G\big)$ are given by the formula
    \begin{equation} \label{T Matrix}
        T_{(K, a, \rho_a), (L, b, \psi_b)} = \delta_{K, L}\delta_{\rho_a, \psi_b}\frac{\chi_{\rho_a}(a)}{\chi_{\rho_a}(e)}.
    \end{equation}
    and
    \begin{align} 
        S_{(K, a, \rho_a), (L, b, \psi_b)} = \frac{1}{|G|}\sum_{\mathclap{\substack{g \in K \\ h \in L \cap C(g)}}} & \left(\frac{\theta_a(a_g, h)\theta_a(a_gh, a_g^{-1})\theta_b(b_h, g)\theta_b(b_hg, b_h^{-1})}{\theta_g(a_g^{-1}, a_g)\theta_h(b_h^{-1}, b_h)}\right)^* \nonumber
        \\ & \quad \times \chi^*_{\rho_a}(a_gha_g^{-1})\chi^*_{\psi_b}(b_hgb_h^{-1}). \label{S Matrix}
    \end{align}   
\end{mythm}

In order to make sense of this theorem and the following discussion, we recall some basic projective character theory. As we will need to work with projective representations controlled by a specified $2$-cocycle rather than merely homomorphisms $G \to \PL(V)$ we provide the (straight forward) proofs in the appendix.

Let $G$ be a finite group, $V$ a complex vector space and $\beta$ a unitary $2$-cocycle.

\begin{mydef} \label{beta proj rep}
    A $\beta$-projective representation is a pair $(V, \rho)$ with $V$ a vector space and $\rho$ a function $\rho: G \to \GL(V)$ satisfying
    \[
        \rho(g)\rho(h) = \beta(g, h)\rho(gh)
    \]
\end{mydef}

Many features from regular representations theory carry over to this slightly more general setting. 

\begin{mydef} \label{proj char}
    Given a projective representation, the projective character $\chi_{\rho}$ of the representation is defined to be
    \[
        \chi_{\rho}(g) = \tr\big(\rho(g)\big).
    \]
    Additionally, for a fixed $\beta$, a $\beta$-projective representation is determined uniquely by its character.
\end{mydef}

Projective characters act similarly to linear characters however many of the classical relations need to be twisted by the $2$-cocycle $\beta$.

\begin{align}
    & \chi_{\rho}(g)\chi_{\rho}(h) = \beta(g, h)\chi_{\rho}(gh) \label{twisted char}
    \\ & \chi_{\rho}(hgh^{-1}) = \frac{\beta(h^{-1}, h)}{\beta(h, g)\beta(hg, h^{-1})} \chi_{\rho}(g) \label{twisted conjugate}
    \\ & \chi_{\rho}(g^{-1}) = \beta(g, g^{-1})\chi^{*}_{\rho}(g), \label{inv char}
\end{align}
where $\chi^{*}$ denotes the complex conjugate of $\chi$.

\subsection{The category of \texorpdfstring{$G$}{G}-graded vector spaces}
    Whilst the categories $\B Z\big(\V^{\omega} G\big)$ are well known, there are two binary choices of conventions that need to be made. These choices, of the direction of the associator morphism and of the half-braiding morphisms, are not entirely settled in the literature. We've resisted the temptation to use the most convenient ones and instead follow \cite{Tens:Cat}. This has the effect that in Theorem \ref{Representation Equivalence Ribbon} there is a perhaps unexpected appearance of $\omega^{-1}$.

    For a group $G$ and unitary three-cocycle $\omega$, the category $\V^{\omega}G$ has objects $G$-graded vector spaces
    \[
         V = \bigoplus_g V_g
    \]
    and morphisms $G$-graded homomorphisms. The simple objects in the category are the one dimensional vectors spaces sitting over each group element
    \[
        (\delta_g)_h = \begin{cases}
            \m{C} & \text{ if } g = h
            \\ 0 & \text{ else}.
        \end{cases}
    \]
    There is a natural tensor product structure given by
    \[
        (V \otimes W)_g = \bigoplus_{g = hk} V_h \otimes W_k
    \]
    with identity object $\delta_e$, trivial unitors and the associator
    \[
        \omega: (\delta_g \otimes \delta_h) \otimes \delta_k \iso \delta_g \otimes (\delta_h \otimes \delta_k)
    \]
    is exactly multiplication by $\omega(g, h, k)$. Note the direction of the associator, this is one of the binary choices we need to make. If you wish to have associator going in the opposite direction you need to use $\omega(g, h, k)^{-1}$. The rigid and pivotal structures on $\V^{\omega}G$ follow naturally.

    Next, we define the Drinfeld center $\B Z\big(\V^{\omega}G\big)$. This has elements pairs $(X, \alpha)$ with $X \in \V^{\omega}G$ and $\alpha$ a half braiding
    \[
        \alpha_Y: Y \otimes X \to X \otimes Y.
    \]
    The morphisms, tensor product, rigid and pivotal structures are directly induced from the relevant structures on $\V^{\omega}G$ and the braiding is given by
    \[
        \sigma_{(X, \alpha), (Y, \beta)} = \beta_X.
    \]
    this is the other binary choice, another common definition is to have the half braidings go in the opposite direction and define the braiding to be $\alpha_Y$. This corresponds to using the inverse braiding $\sigma'_{X, Y} = (\sigma_{Y, X})^{-1}$ and inverts the modular data.

\subsection{An intricate Hopf algebra}
    Let $G$ be a finite group with identity $e$ and $\omega \in H^3(G, \m C)$ a unitary $3$-cocycle. Let
    \begin{equation*} 
        \gamma_{x}(h, l) = \frac{\omega(h, l, x)\omega(x, x^{-1}hx, x^{-1}lx)}{\omega(h, x, x^{-1}lx)}.
    \end{equation*}
    and recall the definition of $\theta_g(x, y)$ given in Equation \ref{Derived two cocycle}.

    Then the ribbon quasi Hopf algebra $D^{\omega}G$ is defined following \cite{Quant:Grp, Altschuler:1992}. Start with the vector space over $\m C$ generated by symbols $\delta_g \ov{x}$ with $g, x \in G$. The algebra structure is

    \begin{gather*}
        \nabla(\delta_{g}\ov{x}, \delta_{h}\ov{y}) = \theta_{g}(x, y)\delta_{g, xhx^{-1}}\delta_{g}\ov{xy} = \begin{cases}
          \theta_{g}(x, y)\delta_{g}\ov{xy} & \text{if } g = xhx^{-1}
          \\ 0 & \text{ otherwise}
        \end{cases} \label{multiplication}
        \\ \eta(k) = k\sum_{g \in G} \delta_g\ov{e}.
    \end{gather*}
    where $\delta_{g, h}$ is the Kronecker delta function. The coalgebra structure is
    \begin{gather*}
        \Delta(\delta_g\ov{x}) = \sum_{h \in G} \gamma_{x}(h, h^{-1}g) \delta_{h}\ov{x} \otimes \delta_{h^{-1}g}\ov{x}
        \\ \epsilon(\delta_g\ov{x}) = \delta_{g, e}.
    \end{gather*}
    This structure is not coassociative but is quasi coassociative with invertible element
    \[
        \Phi = \sum_{g, h, k \in G} \omega(g, h, k)\delta_g\ov{e} \otimes \delta_h\ov{e} \otimes \delta_k\ov{e}.
    \]
    This means that
    \[
        (\Delta \otimes 1) \circ \Delta(h) = \Phi (1 \otimes \Delta) \circ \Delta(h) \Phi^{-1}.
    \]
    Next, the quasi Hopf algebra structures are
    \begin{equation*}
        S(\delta_h\ov{g}) = \theta_{h^{-1}}(g, g^{-1})^{-1} \gamma_g(h, h^{-1})^{-1} \delta_{g^{-1}h^{-1}g}\ov{g^{-1}}
    \end{equation*}
    with $\alpha = \textbf{1} = \eta(1) = \sum_{g \in G} \delta_g\ov{e}$ and $\beta = \sum_{g}\omega(g, g^{-1}, g)\delta_g\ov{e}$. The quasitriangular element is 
    \[
        R = \sum_{g, h \in G}\delta_g \ov{e} \otimes \delta_h \ov{g}.
    \]
    and the ribbon element is
    \[
        v = \sum_{g \in G} \omega(g^{-1}, g, g^{-1}) \delta_g \ov{g^{-1}}.
    \]
    This forms a ribbon quasi Hopf algebra. We omit the simple checks of these conditions; they all follow from the cocycle condition and in general are equivalent to one of the three following identities.
    \begin{gather}
        \theta_g(x, y)\theta_g(xy, z) = \theta_g(x, yz)\theta_{x^{-1}gx}(y, z) \label{2-cocycle:centre}
        \\ \theta_g(x, y)\theta_h(x, y)\gamma_x(g, h)\gamma_y(x^{-1}gx, x^{-1}hx) = \theta_{gh}(x, y)\gamma_{xy}(g, h) \label{theta:gamma:relation}
        \\ \gamma_x(g, h)\gamma_x(gh, k)\omega(x^{-1}gx, x^{-1}hx, x^{-1}kx) = \gamma_x(h, k)\gamma_x(g, hk)\omega(g, h, k) \label{gamma:identity}.
    \end{gather}
    The first identity here, Equation \ref{2-cocycle:centre} shows that for any element $g \in G$, $\theta_g(x, y)$ is a $2$-cocycle on $C(g)$.

    The reason for defining this ribbon quasi Hopf alegbra is due to the following theorem.

    \begin{mythm} \label{Representation Equivalence Ribbon}
        The ribbon fusion categories $\Rep(D^{\omega}G)$ and $\B Z\big(\V^{\omega^{-1}} G\big)$ are equivalent.
    \end{mythm}

    For brevity, the proof of this theorem has been pushed to the appendix.

    \begin{mycor}
        The $S$ and $T$ matrices of $\B Z\big(\V^{\omega^{-1}} G\big)$ and $\Rep(D^{\omega}G)$ are equivalent up to a re-ordering of the simple objects.
    \end{mycor}

    Note that it does not make sense to ask for the $S$ and $T$ matrices to be identical as there is no set ordering of simple objects.

\subsection{Deriving the S and T matrices} \label{derive:S:T}

    Theorem \ref{Representation Equivalence Ribbon} allows us to move between $\B Z\big(\V^{\omega^{-1}} G\big)$ and $\Rep(D^{\omega}G)$. The first step is to classify the simple objects of $\B Z\big(\V^{\omega^{-1}} G\big)$.

\subsubsection{Classifying the Simple Objects} \label{simple:obs}

    The classification of the irreducible left modules of $D^{\omega} G$ can be found in \cite{Altschuler:2004}. Letting $I(G)$ denote the conjugacy classes of $G$, and choosing an element $a \in K$ for each conjugacy class $K \subset I(G)$. The irreducible left modules are classified by triples $(K, a, \rho_a)$ with $\rho_a$  an irreducible $\theta_a$-projective representation of $C(a)$.
    
    This classification comes from analysing the multiplicative structure on $D^{\omega}G$. If $(V, \rho)$ is a representation of $D^{\omega}G$, then for any $g \in G$, $\rho(\delta_g\ov{e})$ must be a projection matrix (since $\nabla(\delta_g\ov{e}, \delta_g\ov{e}) = \delta_g\ov{e}$) and thus diagonalizable. Additionally, for any other $h \in G$,
    \[
        \nabla(\delta_g\ov{e}, \delta_h\ov{e}) = \nabla(\delta_h\ov{e}, \delta_g\ov{e}) = \begin{cases}
        \delta_g \ov{e} & \text{if } g = h
        \\ 0 & \text{ otherwise}.
        \end{cases}
    \]
    Therefore the set of matrices $S = \{\rho(\delta_g\ov{e})\}_{g \in G}$ are simultaneously diagonalisable. Hence, with respect to a basis $\{v_1, \cdots, v_n\}$ diagonalizing $S$, $\rho(\delta_g\ov{e})$ will be a projection onto some subspace $\{v_{g_1}, \cdots, v_{g_{k_g}}\}$. As $\sum_{g \in G} \delta_g\ov{e} = \textbf{1}$, every basis element will be in the image of a unique $\rho(\delta_g\ov{e})$ and so $V$ splits as 
    \[
        V = \bigoplus_{g \in G} V_g
    \]
    where $V_g = \im(\rho(\delta_g\ov{e})) = \Span(v_{g_1}, \cdots, v_{g_{k_g}})$. As
    \[
        \rho(\delta_g\ov{x}) \circ \rho(\delta_{x^{-1}gx}\ov{e}) = \rho(\delta_g \ov{x}) = \rho(\delta_g\ov{e}) \circ \rho(\delta_g\ov{x})
    \]
    we see that $\rho(\delta_g \ov{x})$ is a linear map from $V_{x^{-1}gx} \to V_g$ and acts as $0$ on $V_h$ for $h \neq x^{-1}gx$. As $\delta_g \ov{x}$ is invertible, $V_g$ and $V_{x^{-1}gx}$ also must have the same dimension for all $x \in G$.

    Looking back at Equation \ref{multiplication}, observe that $\nabla(\delta_g\ov{x}, \delta_h\ov{y}) \neq 0$ if and only if $g = xhx^{-1}$. In particular, this means that $g$ and $h$ must be in the same conjugacy class and so the underlying algebra of $D^{\omega}G$ can be written as the direct sum over the conjugacy classes $I(G)$ of $G$
    \[
        D^{\omega}G = \bigoplus_{K \in I(G)} D^{\omega}(K, G).
    \]
    Here $D^{\omega}(K, G)$ is the subalgebra of $D^{\omega}G$ generated by $\delta_g \ov{x}$ for $g \in K$ and $x \in G$ with identity $\sum_{g \in K} \delta_g\ov{e}$. 

    As $D^{\omega}G$ splits into a direct sum of subalgebras, any irreducible left module must entirely sit over one of these $D^{\omega}(K, G)$. Consider again the subalgebra $D^{\omega}(a, C(a))$ inside $D^{\omega}(K, G)$. Observe that this subalgebra is exactly the $\theta_a$-twisted group algebra of $C(a)$ also known as $\m C_{\theta_a}[C(a)]$. Therefore irreducible left modules of $D^{\omega}(a, C(a))$ correspond bijectively to irreducible $\theta_a$-projective representations of $C(a)$.

    Then, if $(V, \rho)$ is a left $D^{\omega}(K, G)$ module, we can define $V_a = \im(\rho(\delta_a \ov{e}))$ which will be $\theta_a$-projective representations of $C(a)$ with action given by $\rho_a = \rho|_{D^{\omega}(a, C(a)) \times V_a}$. This restriction map sends irreducible left modules of $D^{\omega}(K, G)$ to irreducible $\theta_a$-projective representations of $C(a)$ and has an inverse given by an induction functor.

\subsubsection{Deriving the S Matrix}

\begin{figure}
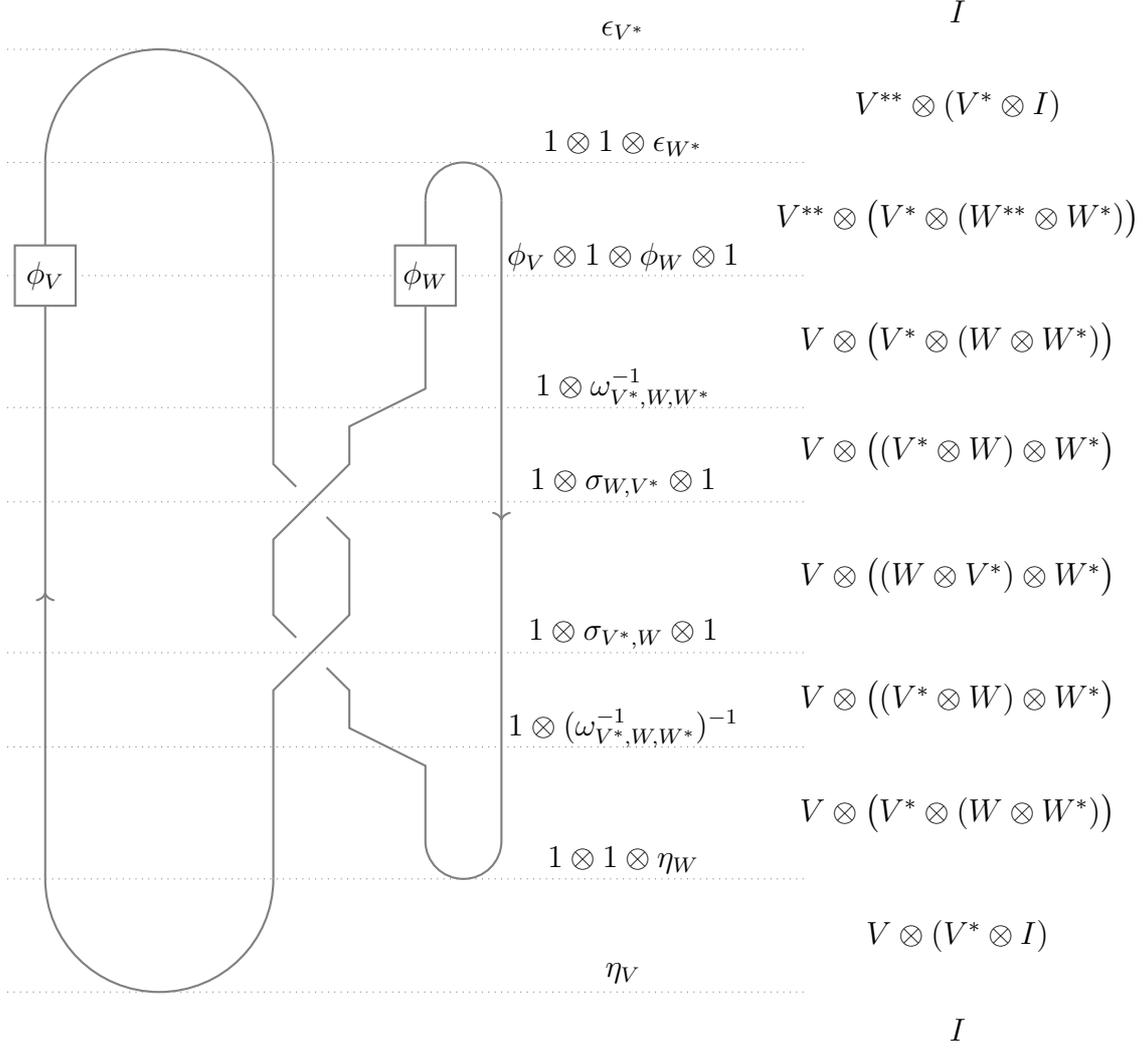
 
        \[
            \ModularDataSDetailed
        \]
        \caption{The string diagram for the entry of the $S$ matrix corresponding to the simple objects $(K, a, \rho_a)$ and $(L, b, \psi_b)$.} \label{S Diagram}
    \end{figure}

    In order to derive the $S$ matrix, we first need to find the corresponding morphism inside $\Rep(D^{\omega}G)(I \to I)$. Let $(K, a, \rho_a)$ and $(L, b, \psi_b)$ be two simple objects in $\B Z(\V^{\omega^{-1}}G)$. These simple objects will be unambiguously referred to by the vector spaces $V$ and $W$ which correspond to the left modules of $D^{\omega}G$ with action given by $\rho$ and $\psi$. Then, the value of $\tilde{S}_{(K, a, \rho_a), (L, b, \psi_b)}$ is the morphism in $\B Z(\V^{\omega^{-1}}G\big)\big(I \to I\big)$ given by Figure \ref{S Diagram}. Note that unitors in Figure \ref{S Diagram} have been omitted because they are trivial. 

    While it is possible to derive the formula for the $S$ matrix directly from Figure \ref{S Diagram}, we can make life much easier for ourselves bu applying some known results. In particular, Altsch\"{u}ler and Coste show in \cite{Altschuler:1992} (On page 98) that the quantum trace in $\Rep(D^{\omega}G)$ is exactly the same as the linear trace. Hence if we can shift everything across to $D^{\omega}G$ we simple need to find the action of the double braiding. We have
    
    For any two left modules $(V, \rho)$ and $(W, \psi)$:
    \[
        \sigma_{V, W} \xrightarrow{\cong} \tau \circ (\rho \otimes \psi)(R) \quad \quad v\otimes w \mapsto \sum_{g, h \in G} \psi(\delta_h \ov{g}, w) \otimes \rho(\delta_g\ov{e}, v)
    \]
    Hence, the double braiding in our diagram corresponds to the map on $V^* \otimes W)$ given by
    \[
        \sigma^2_{V, W} =\tau \circ (\psi \otimes \rho^*)(R) \circ \tau \circ (\rho^* \otimes \psi)(R),
    \]
    which acts on a pair $v^* \otimes w \in V^* \otimes W$ by
    \begin{align*}
        \tau \circ (\psi \otimes \rho^*)(R) \circ \tau \circ (\rho^* \otimes \psi)(R)(v^* \otimes w) & = \tau \circ (\psi \otimes \rho^*)(R)\Big(\sum_{g, h \in G} \psi(\delta_h \ov{g}, w) \otimes \rho^*(\delta_g\ov{e}, v^*)\Big)
        \\ & = \sum_{g, h, k, l \in G} \rho^*\big(\delta_l \ov{k}, \rho^*(\delta_g\ov{e}, v^*)\big) \otimes \psi\big(\delta_k \ov{e}, \psi(\delta_h \ov{g}, w)\big).
    \end{align*}
    As $\psi$ and $\rho^*$ are representations, composition is equivalent to multiplication. Thus, as $\theta_g(x, e) = \theta_g(e, x) = 1$ for all $g, x \in G$,
    \[
        \rho^*(\delta_l \ov{k})\rho^*(\delta_g\ov{e}) = \begin{cases}
            \rho^*(\delta_{kgk^{-1}}\ov{k}) & \text{ if } l = kgk^{-1}
            \\ 0 &  \text{ otherwise}
        \end{cases}
    \]
    and 
    \[
        \psi(\delta_k \ov{e})\psi(\delta_h \ov{g}) = \begin{cases}
            \psi(\delta_h\ov{g}) & \text{ if } k = h
            \\ 0 &  \text{ otherwise}.
        \end{cases}
    \]
    Thus we can remove the sums over $k$ and $l$, replacing $l$ with $hgh^{-1}$ and $k$ with $h$. Additionally, replace $g$ with $h^{-1}gh$. Thus the double braiding acts as 
    \[
        \sigma^2_{V^*, W}(v, w) = \sum_{g, h \in G} \Big(\rho^*(\delta_{g}\ov{h}, v^*), \psi(\delta_h\ov{h^{-1}gh}, w)\Big).
    \]
    Choosing bases for $V$ and $W$ and dual bases for $V^*$ and $W^*$, we apply the definition of the trace to get 
    \begin{align*}
        S_{V, W} = \tr(\sigma^2_{V^*, W}) & = \sum_{\substack{i, j \\ g, h \in G}} (v_i)^{**}\Big(\rho^*(\delta_{g}\ov{h}, v^i)\Big) w^j\Big(\psi(\delta_h\ov{h^{-1}gh}, w_j)\Big)
        \\ & = \sum_{\substack{i, j \\ g, h \in G}} \Big(\rho^*(\delta_{g}\ov{h}, v^i)\Big)v_i \  w^j\Big(\psi(\delta_h\ov{h^{-1}gh}, w_j)\Big)
    \end{align*}
    Now we need to use the definition of the dual representation.
    \[
        \Big(\rho^*(\delta_{g}\ov{h}, v^i)\Big)v_i = v^i\Big(\rho\big(S(\delta_{g}\ov{h}, v_i\big)\Big) = \frac{1}{\theta_{g^{-1}}(h, h^{-1}) \gamma_h(g, g^{-1})} v^i\Big(\rho\big(\delta_{h^{-1}g^{-1}h}\ov{h^{-1}}, v_i\big)\Big).
    \]
    Swapping $g$ with $g^{-1}$ we get
    \[
        S_{V, W} = \sum_{\substack{i, j \\ g, h \in G}} \frac{1}{\theta_{g}(h, h^{-1}) \gamma_h(g^{-1}, g)} v^i\Big(\rho(\delta_{h^{-1}gh}\ov{h^{-1}})(v_i)\Big) w^j\Big(\psi(\delta_h\ov{h^{-1}g^{-1}h})(w_j)\Big)
    \]
    We choose bases for $V$ and $W$ that diagonalise the action of $\delta_g\ov{e}$ for all $g \in G$. Then $w^j\Big(\psi(\delta_h\ov{h^{-1}gh}, w_j)\Big)$ will be zero unless $h^{-1}g^{-1}h \in C(h)$ and $w_j$ lies in $W_h$. This implies that $g \in C(h)$ and that the first factor will be $0$ unless $v_i$ lies in $V_g$. Let $w_{h_j}$ be a basis for $W_h$ and let $v_{g_i}$ be a basis of $V_g$. Then our expression simplifies to
    \[
        S_{V, W} = \sum_{\mathclap{\substack{g \in G \\ h \in C(g) \\ g_i, h_j}}} \frac{1}{\theta_{g}(h, h^{-1})\gamma_{h}(g^{-1}, g)} v^{g_i}\big(\rho(\delta_{g}\ov{h^{-1}}, v_{g_i})\big) w^{h_j}\big(\psi(\delta_{h}\ov{g^{-1}}, w_{h_j})\big).
    \]
    As $V$ and $W$ correspond to irreducible $D^{\omega}G$ modules they are entirely concentrated over conjugacy classes $K$ and $L$. Then $\sum_{g_i} v^{g_i}\big(\rho(\delta_{g}\ov{h^{-1}}, v_{g_i})\big)$ is precisely the trace of the matrix $\rho(\delta_{g}\ov{h^{-1}})$ and, as $\rho(\delta_{g}\ov{h^{-1}})$ acts by $0$ outside $V_g$,
    \[
        \tr\big(\rho(\delta_{g}\ov{h^{-1}})\big) = \tr\big(\rho_g(h^{-1})\big) = \chi_{\rho_g}(h^{-1}).
    \]
    Where $\chi_{\rho_g}$ is a projective character corresponding to the projective representation of $\rho_g = \rho|_{V_g}$. Additionally, as $g$ and $h$ commute, $\gamma_{h}(g^{-1}, g) = \theta_{h}(g^{-1}, g)$ and so
    \[
        S_{V, W} = \sum_{\mathclap{\substack{g \in K \\ h \in KL \cap C(g)}}} \frac{1}{\theta_{g}(h, h^{-1})\theta_{h}(g^{-1}, g)} \chi_{\rho_g}(h^{-1})\chi_{\psi_h}(g^{-1}).
    \]
    Writing these as projective characters allows us to use some projective character theory. In particular we can apply Equation \ref{inv char} which gives
    \[
        \chi_{\rho_g}(h^{-1}) = \theta_{g}(h, h^{-1})\chi^*_{\rho_g}(h).
    \]
    Using this equality and observing that the $2$-cocycle condition for unitary cocycles implies that $\theta_{h}(g^{-1}, g) = \theta_{h}(g, g^{-1})$, we can simplify the expression for the $S$ matrix component further to 
    \[
        \sum_{\mathclap{\substack{g \in K \\ h \in L \cap C(g)}}} \chi^*_{\rho_g}(h)\chi^*_{\psi_h}(g).
    \]
    While this is the simplest version of the formula, it is not useful for computations as we can only compute the representations $\rho_a$ and $\psi_b$. Therefore we need to find a way to express $\chi^*_{\rho_g}(h)$ in terms of a $\rho_a$ character.
    
    Let $y$ be an element conjugate to $a$ and $h \in C(y)$. Then we can pick an element $a_y$ in $G$ witnessing the conjugacy of $a$ and $y$. That is to say, $a = a_yya_y^{-1}$. Then consider the following calculation
    \begin{align}
        \chi_{\rho_y}(h) & = \tr\Big(\rho(\delta_y\ov{h})\Big) \nonumber
        \\ & = \theta_y(a_y^{-1}, a_y)^{-1}\tr\Big(\rho(\delta_y\ov{a_y^{-1}})\rho(\delta_a\ov{a_y})\rho(\delta_y\ov{h})\Big) \nonumber
        \\ & = \theta_y(a_y^{-1}, a_y)^{-1}\tr\Big(\rho(\delta_a\ov{a_y})\rho(\delta_y\ov{h})\rho(\delta_y\ov{a_y^{-1}})\Big) \nonumber
        \\ & = \frac{\theta_a(a_y, h)\theta_a(a_yh, a_y^{-1})}{\theta_y(a_y^{-1}, a_y)}\tr\Big(\rho(\delta_a\ov{a_yha_y^{-1}})\Big) \nonumber
        \\ & = \frac{\theta_a(a_y, h)\theta_a(a_yh, a_y^{-1})}{\theta_y(a_y^{-1}, a_y)}\chi_{\rho_a}(a_yha_y^{-1}). \label{conjugate character}
    \end{align}
    Hence given a fixed $a \in K$ and $b \in L$, for every element $g \in K$ and $h \in L$ choose an $a_g$ and $b_h$ satisfying $a = a_gha_g^{-1}$ and $b = b_hhb_h^{-1}$. Then, applying the formula we just derived to current expression of the $S$ matrix gives us
    \begin{align*}
        S_{(K, a, \rho_a), (L, b, \psi_b)} = \sum_{\mathclap{\substack{g \in K \\ h \in L \cap C(g)}}} & \left(\frac{\theta_a(a_g, h)\theta_a(a_gh, a_g^{-1})\theta_b(b_h, g)\theta_b(b_hg, b_h^{-1})}{\theta_g(a_g^{-1}, a_g)\theta_h(b_h^{-1}, b_h)}\right)^* \chi^*_{\rho_a}(a_gha_g^{-1})\chi^*_{\psi_b}(b_hgb_h^{-1}).
    \end{align*}   
    Dividing through by the normalisation constant, $|G|$, gives the formula in Theorem \ref{S and T matrices}.

\subsubsection{Deriving the T Matrix}

\begin{figure}
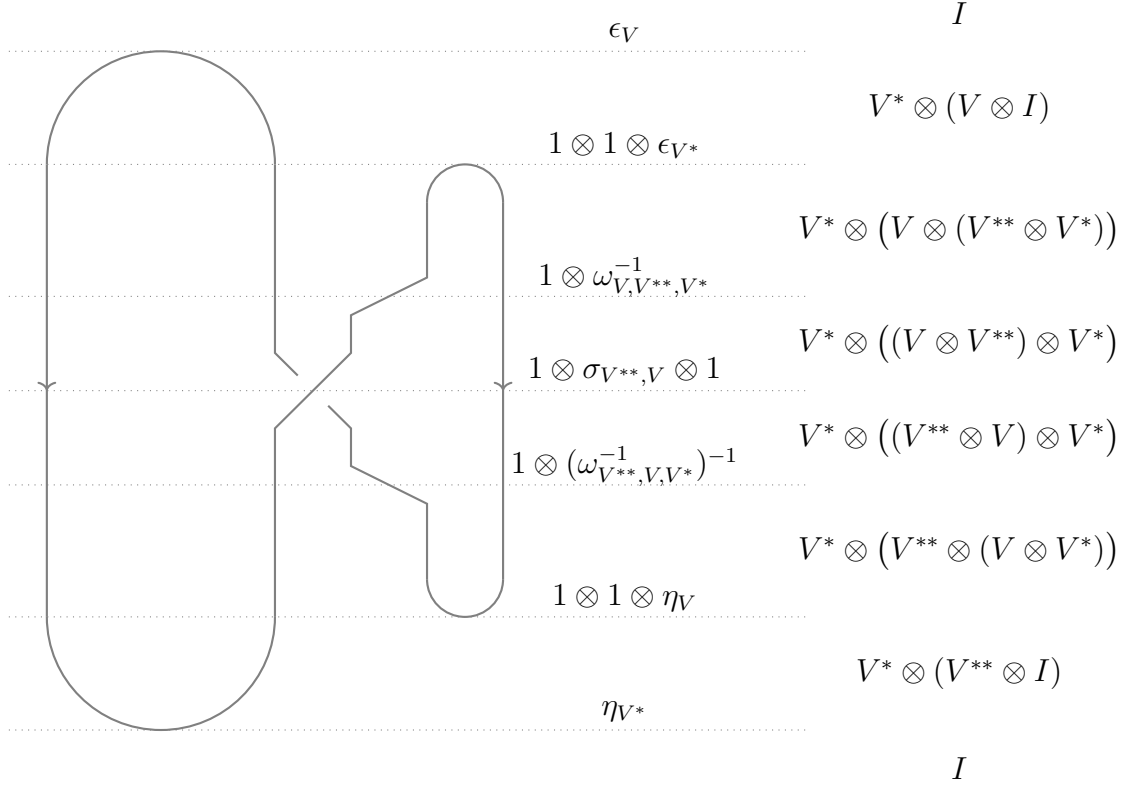
 
        \[
            \ModularDataTDetailed
        \]
        \caption{A detailed string diagram for deriving the diagonal entry of the $T$ matrix corresponding to the simple object $(K, a, \rho_a)$.} \label{T Diagram}
    \end{figure}

    To find the $T$ matrix let us draw another detailed picture as we did for the $S$ matrix. Let $(K, a, \rho_a)$ be a simple object with corresponding left module $(V, \rho)$, then $|V|T_{(K, a, \rho_a)}$ is given by the morphism in Figure \ref{T Diagram}.
    Similarly to our calculation for the $S$ matrix, we don't need to explicitly compute every morphism. The small loop is exactly the ribbon twist and corresponds in $\Rep(D^{\omega} G)$ to the action of the inverse of the ribbon element
    \[
         v^{-1} = \sum_{g \in G} \delta_g \ov{g}.
    \]
    Then, we apply the result from \cite{Altschuler:1992} again which allows us to use the regular trace as opposed to the quantum trace. Hence we find that
    \begin{align*}
        T_V & = \sum_{\substack{i \\ g \in G}} v^i(\rho(\delta_g\ov{g}, v_i)
        \\ & = \sum_{g \in G} \chi_{\rho_g}(g)
    \end{align*}
    Next, recall the chosen element $a \in K$ and pick an $a_g \in G$ such that $a = a_gga_g^{-1}$. Then treating $\sum_{g_i} \Big(\rho(\delta_{g}\ov{g}, v_{g_i})\Big)$ as a character and using Equation \ref{conjugate character} again to replace $\chi_{\rho_g}$ with a $\chi_{\rho_a}$ we find that
    \[
        T_{(K, a, \rho_a), (L, b, \psi_b)} = \frac{\delta_{K, L}\delta_{\rho_a, \psi_b}}{|V|} \sum_{g\in K} \chi_{\rho_g}(g) = \frac{\delta_{K, L}\delta_{\rho_a, \psi_b}}{|V|} \chi_{\rho_a}(a) \sum_{g\in K} \frac{\theta_a(a_g, g)\theta_a(a_g g, a_g^{-1})}{\theta_g(a_g^{-1}, a_g)}.
    \]
    Consider further the term
    \[
        \frac{\theta_a(a_g, g)\theta_a(a_gg, a_g^{-1})}{\theta_g(a_g^{-1}, a_g)}.
    \]
    First observe that Equation \ref{2-cocycle:centre} with $g = g$, $x = a^{-1}_g$, $y = a^{-1}_g$, and $z = g$ gives the equality
    \[
        \theta_g(a^{-1}_g, a_g) \theta_g(e, g) = \theta_g(a^{-1}_g, a_gg)\theta_{a_gga^{-1}_g}(a_g, g).
    \]
    As $\theta$ is normalised, and $a_gga^{-1}_g = a$, this summation term can be rearranged by
    \begin{align*}
        \frac{\theta_a(a_g, g)\theta_a(a_gg, a_g^{-1})}{\theta_g(a_g^{-1}, a_g)} & = \frac{\theta_a(a_gg, a_g^{-1})}{\theta_g(a^{-1}_g, a_gg)}
        \\ & = \frac{\omega(a, a_g,g, a_g^{-1})\omega(a_gg, a_g^{-1}, a_gg^{-1}a_g^{-1}aa_gga_g^{-1})\omega(a_g^{-1}, a_gga^{-1}_g, a_gg)}{\omega(a_gg, g^{-1}a_g^{-1}aa_gg, a_g^{-1})\omega(g, a^{-1}_g, a_gg)\omega(a_g^{-1}, a_gg, g^{-1}a_g^{-1}a_gga^{-1}_ga_gg)}
        \\ & = \frac{\omega(a_gga_g^{-1}, a_g,g, a_g^{-1})\omega(a_gg, a_g^{-1}, a_gga_g^{-1})\omega(a_g^{-1}, a_gga^{-1}_g, a_gg)}{\omega(a_gg, g, a_g^{-1})\omega(g, a^{-1}_g, a_gg)\omega(a_g^{-1}, a_gg, g)}
    \end{align*}
    This expression appears complicated but it can be simplified using the $3$-cocycle condition. Applying the $3$-cocycle condition on the tuple $(a_gg, a^{-1}_g, a_gg, a^{-1}_g)$ simplifies this expression to
    \[
        \frac{\omega(a_gg, a_g^{-1}, a_g,g)\omega(a_g^{-1}, a_gg, a_g^{-1})\omega(a_g^{-1}, a_gga^{-1}_g, a_gg)}{\omega(g, a^{-1}_g, a_gg)\omega(a_g^{-1}, a_gg, g)}
    \]
    which the $3$-cocycle condition on the tuple $(a^{-1}_g, a_gg, a^{-1}_g, a_gg)$ shows must be $1$. Therefore, as $|V| = |K||V_a| = |K|\chi_{\rho_a}(e)$ we are left with the simple formula
    \begin{equation}
        T_{(K, a, \rho_a), (L, b, \psi_b)} = \delta_{K, L}\delta_{\rho_a, \psi_b} \frac{\chi_{\rho_a}(a)}{\chi_{\rho_a}(e)}.
    \end{equation}
    This completes the proof of Theorem \ref{S and T matrices}.
\section{Computing Modular Data in GAP}
\label{cha:CompData}

Now that we have these formula for the S and T matrices, we can proceed to writing code to produce the modular data. The programming language GAP \cite{GAP4}, was chosen for this project due to both its speed and inbuilt functions. In particular, GAP can easily create groups from a given presentation and find the corresponding character tables. On top of GAP there were two packages used, `hap' and `IO'. The package `hap' \cite{HAP}, adds in numerous functions relating to group cohomology and was required for the code written by Mignard and Schauenburg. The package `IO' \cite{IO} was used to aid the back end of storing the computed modular data.

The code discussed here is available with the arXiv sources or alternatively at \url{https://tqft.net/web/research/students/AngusGruen}. A database of computed modular data for all groups with order less than 48 can also be found at this website.

Given a group $G$ and a cocycle $\omega$, there are two steps that need to be completed in order to find the modular data. The first step is to create a list of all the simple objects of $\B Z\big(\V^{\omega^{-1}}G\big)$ and the second is to create the $S$ and $T$ matrices from this list.

\subsection{Computing the Simple Objects}

   The simple objects of $\B Z\big(\V^{\omega^{-1}}G\big)$ correspond to tuples $(K, a, \rho_a)$ with $K$ a conjugacy class of $G$ with representative $a$ and $\rho_a$ an irreducible $\theta_a$-projective representation of $C(a)$.

  The function \texttt{GenerateSimpleObjects} begins by producing the set of conjugacy classes of $G$. Then, for each conjugacy class $K$, it will pick a representative $a$ and create a list of pairs\footnote{Recall that $a_g$ is chosen so that $a = a_gga_g^{-1}$.} $[g, a_g]$ for every element in $K$. Next we need to find irreducible projective representations. In order to accomplish this we use the following theorems about projective representations. A proof of these results can be found in Appendix \ref{app:app1}.

  \begin{mythm} \label{Group Extension Theorem}
    Let $G$ be a group and $\beta$ a unitary $2$-cocycle. Then there exists a group $G_{\beta}$ called the group extension of $G$ and an injective map $f$ from $\beta$-projective representations of $G$ to linear representations of $G_{\beta}$.
  \end{mythm}

  A more natural construction of $G_{\beta}$ can be found in Appendix \ref{app:app1} but for our purposes here we present a construction which produces a finite presentation $G_{\beta}$ from a finite presentation of $G$. This was given by Flannery and O'Brien in \cite{flannery:2000}.

  As $\beta$ is unitary we can consider it mapping into integers mod $m$ for some $m$.

  Let $\langle g_i | r_j\rangle$ be a finite presentation of $G$. Then define $\ov{\beta}$ to be a map from words over the alphabet $\{g_i, g_i^{-1}\}$ to $\{0, \cdots, m - 1\}$ by
  \[
      \ov{\beta}(g_{i_1}^{e_1}\cdots g_{i_n}^{e_n}) = \prod_{j = 1}^{n} \beta(g_{i_1}^{e_1}\cdots g_{i_{j - 1}}^{e_{j-1}}, g_{i_j}^{e_j}) \beta(g_{i_j}, g_{i_j}^{-1})^{\frac{e_i - 1}{2}}.
  \]
  A finite presentation for $G_{\beta}$ is then given by
  \[
      G_{\beta} = \Big\langle g_i, x | x^{m}, xg_ix^{-1}g_i^{-1}, r_jx^{-\ov{\beta}(r_j)} \Big\rangle.
  \]

  Using this definition of $G_{\beta}$, given a projective representation $\rho$, $f(\rho)$ is the linear representation $G_{\beta}$ defined by
  \[
    f(\rho)(g_i) = \rho(g_i) \text{ and } f(\rho)(x) = e^{\frac{2\pi i}{m}}.
  \]

  This map $f$ has a very well behaved inverse on $\im(f)$ and it is easy to check if an irreducible representation of $G_{\beta}$ is in $\im(f)$. Specifically given a linear representation $(\psi, V)$ of $G_{\beta}$, $\psi$ is in $\im(f)$ if and only if $\psi(x) = e^{\frac{2\pi i}{m}}$ and if this is the case then $f^{-1}(\psi)$'s is the $\beta$-projective representation of $G$ defined by
  \[
    f^{-1}(\psi)(g_i) = \psi(g_i)
  \]

  We also have the following useful theorem about how $f$ acts on irreducible representations.
  
  \begin{mycor} \label{irr:irr:correspondance}
      The $\beta$-projective representation $(V, \rho)$ will be irreducible if and only if $(V, f(\rho))$ is an irreducible linear representation.
  \end{mycor}

  Picking a conjugacy class $K$ with chosen element $a$, we construct the group extension $C(a)_{\theta_a}$ using the presentation above. Then, the function \texttt{ProjectiveCharacters} finds all $\theta_a$-irreducible characters of $C(a)$ using Corollary \ref{irr:irr:correspondance} and the discussion  above it.

  Next we assemble the list of simple objects, each given as a pair of a conjugacy class and an irreducible projective representation. For each simple object we record $4$ pieces of data, the representative element of the conjugacy class, the projective character of the projective representation, the list of elements $[g, a_g]$ mentioned earlier and the dimension of the projective representation. For speed purposes, the projective character is stored as a lookup dictionary linking each element of the centralizer to the projective character of that element.

  In constructing this list we treat the conjugacy class $\{e\}$ slightly differently for two reasons. Firstly as $\theta_e = 1$, the projective representations correspond to ordinary representations of $G$ which can be computed more quickly and easily. Additionally, we can then guarantee that the first simple object will be the unit object and so the first row and column of the $S$ matrix will be strictly positive rational numbers and correspond to a list of the normalised dimensions of each simple object\footnote{Recall that the normalised dimension of a simple object $V$ in a category $\C C$ is $\frac{\dim(V)}{\sqrt{\dim(\C C)}}$}.

  We also treat the cyclic case $G \cong C_n$ differently. In this case for every $m \in C_n$, all $\theta_m$-irreducible representations of $C(m) = C_n$ are one dimensional\footnote{Note that, unlike in classical representation theory, non-cyclic abelian groups may have $\theta_m$-irreducible representations which are not one dimensional and so this technique does not extend to all abelian groups. For an example of this, look at $G = (\m{Z}/2\m{Z})^5$.}. This observation allows us to directly write down all possible characters of a generator $1 \in C_n$ as
  \[
  		\rho(1)^n = \left(\prod_{i=1}^{n-1} \theta_m(1, i)\right)\rho(n) = \prod_{i=1}^{n-1} \theta_m(1, i)
  \]
  and the $n$ solutions for this equation each correspond to their own irreducible representation.

  It is important to note that some of the inbuilt GAP functions used here are non-deterministic. This means that running this piece of code twice may produce different permutations of the list of simple objects. In particular this means that the $S$ and $T$ matrices produced will not necessarily by identical. They will of course be equivalent up to conjugative by a permutation matrix.

  After generating the simple objects, the code moves on to creating the $S$ and $T$ matrices using the function \texttt{GenerateModularData}.

\subsection{Computing the S and T Matrices}
  
  The simplest approach to computing the $S$ matrix is to directly implement the formula from Theorem \ref{S and T matrices}. The function \texttt{CoefficientSSum} is exactly calculating the coefficient
  \[
    \frac{\theta_a(a_g, h)\theta_a(a_gh, a_g^{-1})\theta_b(b_h, g)\theta_b(b_hg, b_h^{-1})}{\theta_g(a_g^{-1}, a_g)\theta_h(b_h^{-1}, b_h)}.
  \]
  Then, \texttt{sVal} implements Equation \ref{S Matrix} when given $2$ simple objects. To make the $S$ matrix, one could then iterate over the created list of simple objects twice and call \texttt{sVal} on each pair to create the corresponding entry of the matrix. The $T$ matrix is comparatively easier to make as we only need to compute the diagonal. Therefore we only need to loop over the simple objects once, directly implementing Equation \ref{T Matrix}.

  In practice, this code can be sped up a lot. This comes from the fact that calling \texttt{sVal} is expensive and the S matrix has several symmetries.

  The first trick is observing that when $G$ is abelian, all conjugacy classes contain only $1$ element and all centralizers are the entire group. This simplifies the formula for elements of the $S$ matrix to be
  \[
    S_{(\{g\}, g, \rho_g), (\{h\}, h, \psi_h)} = \chi^*_{\rho_g}(h)\chi^*_{\psi_h}(g).
  \]
  Additionally, as $S = S^T$ we only need to compute the upper triangular half of $S$.

  When $G$ is not abelian, more elaborate methods are needed. We use the following result found in \cite{COSTE1994316},
  \begin{mylem} \label{galois:sym}
    Given any row $S_x$ of the $S$ matrix, and any element $\sigma \in \text{Gal}(\m Q[S_x]/\m Q)$, denote $\sigma(S_x)$, the list given by applying $\sigma$ to every element of $S_x$. Then either $\sigma(S_x)$ or $-\sigma(S_x)$ also appears as a row in the $S$ matrix. 
  \end{mylem}
  For our specific case, as $S_1$ is always a positive rational, $-\sigma(S_x)$ will never be a row and so $\sigma(S_x)$ will always be one. Additionally, $\text{Gal}(\m Q[S_x]/\m Q) \subset \text{Gal}(\m Q[T]/\m Q)$ and so this group can be precomputed after calculating the diagonal of the $T$-matrix.

  This allows for the following algorithm. Start with an empty $S$-matrix. Then repeat the following steps:
  \begin{enumerate}
    \item Compute a random unknown row and add it to the $S$ matrix. 
    \item If this computed row is not a Galois conjugate of an already computed row, save all of its Galois conjugates to a list containing other rows that are Galois conjugates of computed rows.
    \item Go through the list of saved Galois conjugates and attempt to place them into the partially filled $S$ matrix using the constraint that $S = S^T$. If there is only a single place that a row could go, place it there.
    \item If the $S$ matrix is not yet complete, return to step $1$.
  \end{enumerate}

  The reason for computing a random row in step 1 as opposed to going through the list of rows deterministically is because we want to avoid computing rows that are equal to Galois conjugates of rows already computed. In practice it is usually the case that many of the Galois conjugates of a row are grouped together and so it is a lot faster in general to compute a random row each time.

  In practice this algorithm can be hundreds of times faster than the direct approach. Exactly what the speed differential is depends on the specific group and $3$-cocycle but even for small groups is it usually at least twice as fast.

\subsection{Constructing the Database}

  Using this code to compute, for a group $G$ and a three cocycle $\omega$, the modular data of $\B Z\big(\V^{\omega^{-1}}G\big)$ we can begin constructing a database of all possible sets of modular data for groups and cocycles below some order. A problem we run into immediately is that for many groups, $|H^3(G)|$ can be enormous. Luckily the following lemma will allow us to massively cut down on the number of different $3$-cocycles we need to consider for each group.

  \begin{mylem}
      If $\omega \in H^3(G)$, $\theta \in \Aut(G)$ then $\V^{\omega}G \overset{\otimes}{\cong} \V^{\theta^*(\omega)}G$.
  \end{mylem}

  Here, $\theta^*(\omega)$ denotes the pullback of $\omega$ by $\theta$. This lemma is particularly helpful when dealing with groups such as powers of $\m Z/2\m Z$. While $\left|H^{3}(\left(\m Z/2\m Z\right)^5, \m{C}^{\times})\right| = 2^{25} = 33554432$, the number of orbits under the action of the automorphism group is merely $88$. This means we simply need to write code to produce a unitary cocycle representatives of $H^3(G, \m C)/\Aut(G)$.

  Symmetry arguments show that if $\{\omega_i\}$ is a complete set of representatives of $H^3(G, \m C)/\Aut(G)$ then $\{\omega^{-1}_i\}$ is also a complete set of representatives. Thus the constructed database will contain modular data for every equivalence class of categories of the form $Z(\V^{\omega} G)$ with $|G|$ less than $48$.

  The code used the algorithm written by Mignard and Schauenburg in their recent paper \cite{mignard:2017}. This computes $H^3(G, \m C)$ using the Universal Coefficient Theorem, produces a basis for the cohomology space and then represents each unique unitary cocycle by a vector. Their code also contains a function to find the orbits of the cohomology space under the action of $\Aut(G)$. For most groups this worked fine but there a couple of groups for which their code was too slow. For groups of order $2^n$, we wrote some bit-flipping code that computed the orbits more quickly. Currently we have not expanded this to groups whose orders are not powers of $2$ but this would be an important exercise if we wanted to increase the size of the modular data database to $48$ and beyond. (In particular to deal with the group $(\m{Z}/2\m{Z})^4\m{Z}/3\m{Z}$)

  This code returns a map $\omega$ from $G \times G \times G \to \{0, \cdots, exp - 1|\} \subset \m Z$ where $exp$ is the LCM of the torsion coefficients of $H^3(G, \m C)$ and $g$, the GCD of elements in the image of $\omega$, identifying $0$ with exp. As all elements in the image of $\omega$ will be divisible by $g$, when viewed as a map into $\m C^{\times}$,
  \[
    \frac{\omega}{g} : G \times G \times G \to \{0, \cdots, \frac{exp}{g} - 1\}
  \] 
  is the same map as $\omega$ but, due to its smaller image, this second map will be quicker to work with. These generated $G$ and $\omega$ can be directly plugged into to the function \texttt{ComputeModularData} which will compute the $S$ and $T$ matrices of $\B Z(\V^{\omega^{-1}})$.
\section{Analysis of Modular Data}
\label{cha:ana}

    While the main aim of this project was to construct the database and thus allow other people to do more detailed analysis, we present some preliminary results here. One interesting question relates to the ranks of $\B Z(\V^{\omega} G)$ given the order of $G$. Recall that the rank of $\B Z(\V^{\omega} G)$ is exactly the dimension of the $S$ and $T$ matrices. Interestingly, for the data that we have constructed, the ranks appear to be heavily restricted by the order of $G$. As $\dim\big(\B Z(\V^{\omega} G)\big) = |G|^2$ we know that $\rank\big(\B Z(\V^{\omega} G)\big)$ is bounded by $|G|^2$ but, in general, only few of the values in the range $\{1, \cdots, |G|^2\}$ are obtained. Figure \ref{table of cats} gives the possible ranks of $\B Z(\V^{\omega} G)$ for a fixed value of $|G|$ it also gives the multiplicities of these ranks with respect to the equivalence classes of $\V^{\omega} G$, the equivalences classes of the modular data of $\B Z(\V^{\omega} G)$ and the equivalence classes obtained by only considering the $T$ matrix of $\B Z(\V^{\omega} G)$.

    Recall that two sets of modular data $S, T$ and $S', T'$ are equivalent if exists a permutation matrix $P$ such that $S = PS'P^{-1}$ and $T = PT'P^{-1}$. If we only considering the $T$ matrix, this is relatively easy to accomplish. Simply sort both diagonals and equality unfortunately this approach is not feasible to fins $S$ matrix equivalences. We wrote to convert this problem to a graph isomorphism question and then used the program \texttt{nauty} written by Brendan Mackay and Adolfo Piperno \cite{McKay201494}. This program \texttt{nauty} accepts a pair of graphs and returns either an automorphism or states that no automorphism exists. The conversion of a matrix to an edge coloured graph is explained on page 60 of the user manual \cite{Nauty}.
    \begin{figure}
        \centering
        \scriptsize
        \begin{tabular} {||p{5mm}|p{32mm}|p{32mm}|p{32mm}|p{32mm}||}
            \hline
            $|G|$ & Possible Ranks of the center & Number of possibly distinct centers of each rank & Number of distinct centers up to modular data equivalence & Number of distinct centers up to $T$ Matrix equivalence \\
            \hline \hline
            2 & 4 & 2 & 2 & 2                           \\ \hline
            3 & 9 & 3 & 3 & 3                           \\ \hline
            4 & 16 & 8 & 7 & 7                          \\ \hline
            5 & 25 & 3 & 3 & 3                          \\ \hline
            6 & 8, 36 & 6, 6 & 6, 6 & 6, 6              \\ \hline
            7 & 49 & 3 & 3 & 3                          \\ \hline
            8 & 22, 64 & 25, 22 & 20, 18 & 16, 17       \\ \hline
            9 & 81 & 10 & 9 & 9                         \\ \hline
            10 & 16, 100 & 6, 6 & 6, 6 & 6, 6           \\ \hline
            11 & 121 & 3 & 3 & 3                        \\ \hline
            12 & 14, 32, 144 & 6, 30, 24 & 6, 27, 21 & 6, 27, 21            \\ \hline
            13 & 169 & 3 & 3 & 3                        \\ \hline
            14 & 28, 196 & 6, 6 & 6, 6 & 6, 6           \\ \hline
            15 & 225 & 9 & 9 & 9                        \\ \hline
            16 & 46, 88, 256 & 66, 189, 73 & 58, 125, 47 & 50, 106, 44      \\ \hline
            17 & 289 & 3 & 3 & 3                        \\ \hline
            18 & 44, 72, 324 & 20, 18, 20 & 18, 18, 18 & 18, 18, 18         \\ \hline
            19 & 361 & 3 & 3 & 3                        \\ \hline
            20 & 22, 64, 400 & 4, 30, 24 & 4, 27, 21 & 4, 27, 21            \\ \hline
            21 & 25, 441 & 3, 9 & 3, 9 & 3, 9           \\ \hline
            22 & 64, 484 & 6, 6 & 6, 6 & 6, 6           \\ \hline
            23 & 529 & 3 & 3 & 3                        \\ \hline
            24 & 21, 42, 56, 86, 128, 198, 576 & 24, 36, 12, 141, 120, 75, 66 & 24, 30, 12, 114, 99, 60, 54 & 24, 18, 9, 102, 96, 48, 51     \\ \hline
            25 & 625 & 10 & 9 & 9                       \\ \hline
            26 & 88, 676 & 6, 6 & 6, 6 & 6, 6           \\ \hline
            27 & 105, 729 & 31, 30 & 23, 24 & 23, 24    \\ \hline
            28 & 112, 784 & 30, 24 & 27, 21 & 27, 21    \\ \hline
            29 & 841 & 3 & 3 & 3                        \\ \hline
            30 & 116, 144, 200, 900 & 18, 18, 18, 18 & 18, 18, 18, 18 & 18, 18, 18, 18 \\ \hline
            31 & 961 & 3 & 3 & 3                        \\ \hline
            32 & 79, 100, 121, 142, 184, 352, 1024 & 60, 589, 72, 129, 1978, 1081, 172 & 60, 330, 72, 94, 1072, 580, 108 & 60, 201, 52, 83, 672, 471, 99     \\ \hline
            33 & 1089 & 9 & 9 & 9                       \\ \hline
            34 & 148, 1156 & 6, 6 & 6, 6 & 6, 6         \\ \hline
            35 & 1225 & 9 & 9 & 9                       \\ \hline
            36 & 36, 64, 126, 176, 288, 1296 & 8, 42, 30, 100, 90, 80 & 8, 42, 28, 81, 81, 63 & 8, 42, 28, 81, 81, 63                       \\ \hline
            37 & 1369 & 3 & 3 & 3                       \\ \hline
            38 & 184, 1444 & 6, 6 & 6, 6 & 6, 6         \\ \hline
            39 & 65, 1521 & 3, 9 & 3, 9 & 3, 9          \\ \hline
            40 & 88, 214, 256, 550, 1600 & 20, 141, 120, 75, 66 & 18, 114, 99, 60, 54 & 18, 102, 96, 48, 51                             \\ \hline
            41 & 1681 & 3 & 3 & 3                       \\ \hline
            42 & 44, 100, 224, 252, 392, 1764 & 6, 6, 18, 18, 18, 18 & 6, 6, 18, 18, 18, 18 & 6, 6, 18, 18, 18, 18                        \\ \hline
            43 & 1849 & 3 & 3 & 3                       \\ \hline
            44 & 256, 1936 & 30, 24 & 27, 21 & 27, 21   \\ \hline
            45 & 2025 & 30 & 27 & 27                    \\ \hline
            46 & 268, 2116 & 6, 6 & 6, 6 & 6, 6         \\ \hline
            47 & 2209 & 3 & 3 & 3                       \\ \hline
        \end{tabular}
        \caption{The distribution of the ranks of $\B Z(\V^{\omega} G)$ for a fixed order of $G$ between $2$ and $47$. (This table can be extended up $|G| = 63$)} \label{table of cats}
    \end{figure}

    There are some clear patterns that can be seen in these results. For example, for all odd primes $p$, there is only $1$ group of order $p$ which has exactly $3$ orbits of $H^3(G, \m C)/\Aut(G)$ each of which corresponds to a different Morita equivalence class. Note that we know that each orbits corresponds to a different Morita equivalence class because the numbers in columns $3$ and $4$ are identical. Additionally, when $p$ is 2 times an odd prime, there are two groups of order $p$ both with $6$ orbits of $H^3(G, \m C)/\Aut(G)$ with each orbit corresponding to a different Morita equivalence class. There also appears to be a a pattern when $|G|$ is three times a prime larger than $3$. 

    As for some more general observations, at every order $|G|$ numerous categories have rank $|G|^2$. This is to be expected as for any abelian group $G$, $\rank\big(\B Z(\V \ G)\big) = |G|^2$. Interestingly though, for the majority of groups calculated, the rank of the centre is independent of the cocycle chosen. This may be an artifact of small groups however as, as $|G|$ increases, the percentage of groups exhibiting at least $2$ different ranks of their centre appears to rise. 

    Another interesting feature is that, for all groups so far computed, $\rank\big(\B Z(\V \ G)\big)$ is an upper bound to $\rank\big(\B Z(\V^{\omega} G)\big)$. It is clear that this should be the case for abelian groups, but it is far less clear if this should be expected to hold for non-abelian groups.

    The third and fourth columns show some interesting characteristics of the equivalence classes of the $S$ and $T$ matrices. In particular they show that for our dataset, when $|G|$ is not divisible by $8$, the equivalence classes of just the $T$ matrix and both the $S$ and $T$ matrices are identical. This is particularly surprising as it remains true at some larger orders of $|G|$ such as $36$ where there are $303$ equivalent classes of both $S$ and $T$ matrices and just the $T$ matrix.

    We can sum together the multiplicities of different ranks of $\B Z(\V^{\omega} G)$ to find upper and lower bounds of the number of Morita equivalent classes of $Z(\V^{\omega} G)$ at different orders of $|G|$. This is shown in Table \ref{Cats Up to Morita}. The naive upper bounds originate directly from the number of equivalence classes of $\V^{\omega} G$ and the lower bounds are the number of inequivalent pairs of modular data.
    \begin{figure}
        \centering
        \begin{tabular} {||c|c|c|c||}
            \hline
            $|G|$ & \# of Groups & Upper bound of \# of $\B Z(\V^{\omega}G)$ & Lower Bound of \# of $\B Z(\V^{\omega}G)$ \\
            \hline \hline
            2 & 1 & 2 & 2 \\ \hline
            $p$ & 1 & 3 & 3 \\ \hline
            $2 \times p$ & 2 & 12 & 12 \\ \hline
            4 & 2 & 8 & 7 \\ \hline
            8 & 5 & 47 & 38 \\ \hline
            9 & 2 & 10 & 9 \\ \hline
            12 & 5 & 60 & 54 \\ \hline
            15 & 1 & 9 & 9 \\ \hline
            16 & 14 & 328 & 230 \\ \hline
            18 & 5 & 58 & 54 \\ \hline
            20 & 5 & 58 & 52 \\ \hline
            21 & 2 & 12 & 12 \\ \hline
            24 & 15 & 474 & 393 \\ \hline
            25 & 2 & 10 & 9 \\ \hline
            27 & 5 & 61 & 47 \\ \hline
            28 & 4 & 54 & 48 \\ \hline
            30 & 4 & 72 & 72 \\ \hline
            32 & 51 & 4081 & 2316 \\ \hline
            33 & 1 & 9 & 9 \\ \hline
            35 & 1 & 9 & 9 \\ \hline
            36 & 14 & 350 & 303 \\ \hline
            39 & 2 & 12 & 12 \\ \hline
            40 & 14 & 422 & 345 \\ \hline
            42 & 6 & 84 & 84\\ \hline
            44 & 6 & 54 & 48\\ \hline
            45 & 1 & 30 & 27 \\ \hline
            48 & 52 & 4422 & 3207 \\ \hline
            49 & 2 & 10 & 9 \\ \hline
            50 & 5 & 58 & 54 \\ \hline
            51 & 1 & 9 & 9 \\ \hline
            52 & 5 & 58 & 52 \\ \hline
            54 & 15 & 434 & 348 \\ \hline
            55 & 2 & 14 & 12 \\ \hline
            56 & 13 & 408 & 333 \\ \hline
            57 & 2 & 12 & 12 \\ \hline
            60 & 13 & 474 & 438 \\ \hline
            63 & 4 & 45 & 41 \\ \hline
        \end{tabular}
        \caption{The lower and upper bounds for the number of Morita equivalence classes of pointed fusion categories of ranks $2$ through $63$. In this table $p$ refers to any odd prime as there are some trends we observed in our data. Additionally note that for $|G| = 55$, the upper bound is known to be correct \cite{Schauenburg:2017}.} \label{Cats Up to Morita}
    \end{figure}

    For $|G|$ between $2$ and $31$ this table is identical to a table in Mignard and Schauenburg's paper \cite{mignard:2017} and in this paper they prove that for $|G| \leq 31$, the lower bound is exact. At $|G| = 32$ we have improved on the lower bound given by Mignard and Schauenburg from $2315$ to $2316$ and for $33 \leq |G| \leq 63$ these lower bounds had not previously been obtained. Due to the exactness of the lower bounds for $|G| \leq 31$, it is likely that the bounds for $|G| \geq 32$ are either optimal or very close to optimal. 

    The invariants used by Mignard and Schauenburg were the Frobenius-Schur indicators and the $T$ matrix. Therefore the improvement in the lower bound as $|G| = 32$ means that the $S$ and $T$ matrices are strictly stronger invariants than the Frobenius-Schur indicators and the $T$ matrix. Unfortunately, as we do not have the equivalence files produced by Mignard and Schauenburg while we know that an example exists, we cannot yet explicitly give a pair of categories which show this fact. It would be relatively simple to compute the Frobenius-Schur indicators from our data and thus pin down this example exactly but this is beyond the scope of this paper. While the result that the $S$ and $T$ matrices are strictly stronger invariants than the Frobenius-Schur indicators and the $T$ matrix was already known \cite{keilberg:2017}, the counterexample given in \cite{keilberg:2017} involves more exotic categories than the ones presented here and so it would be interesting to explicitly find the counterexample at $|G| = 32$.

    There were two effects that prevented us from calculating all of the modular data for groups of order bigger than $63$. The first issue was related to groups with large automorphism groups. While computing the cohomology classes remains quick, finding the orbits of $H^3(G, \m C)/\Aut(G)$ becomes very slow. In particular at order $64$, we have groups such as $\big(\m Z/2\m Z\big)^6, \big(\m Z/2\m Z\big)^4\times\big(\m Z/4\m Z\big)$ as well as non abelian groups with similar structures. This could probably be fixed by a lower level approach, for example the $\big(\m Z/2\m Z\big)^6$ case can be solved with some bit flipping code but we did not pursue this. The other problem is simply the litany of groups at order $64$. There are 319 finite groups with $1 \leq |G| \leq 63$ and $267$ finite groups with $|G| = 64$. Given the number of computing hours it took to extend the database merely from 47 to 63, extending past 64 seems unnecessary at this point.

    If we avoid abelian groups and orders highly divisible by 2 (e.g. 64, 80, 96, 112, 128) the code is able to compute the $S$ and $T$ matrices for groups with order much greater than $63$. For example, the modular data for the groups $S_5, \m{Z}/5\m{Z} \rtimes \m{Z}/31\m{Z}$ and $\PSL(3, 2)$ which have orders $120, 155$ and $168$ respectively are all available in our online database. This middle group $\m{Z}/5\m{Z} \rtimes \m{Z}/31\m{Z}$ is particularly interesting as it is the second smallest example (after $\m{Z}/5\m{Z} \rtimes \m{Z}/11\m{Z}$) where we expect to find (at least) two $3$-cocycles $\omega, \omega'$ so $\B{Z}(\V^{\omega} G) \ncong \B{Z}(\V^{\omega'} G)$, but nevertheless they have equivalent modular data \cite{Schauenburg:2017}.

\appendix

\section{Projective Representation Theory} \label{app:app1}

    Much of this appendix can be found in the paper ``A Character Theory for Projective Representations of Finite Groups'' by Cheng \cite{cheng:1996}.

    Let $G$ be a finite group, let $V$ be a complex vector space, $\beta$ be a unitary $2$-cocycle and $\rho$ be a $\beta$-projective representation of $G$. Recall that this means that $\rho$ is a map $G \to \GL(V)$ satisfying 
     \[
        \rho(g)\rho(h) = \beta(g, h)\rho(gh).
    \]
    
    As $\beta$ is a $2$-cocycle $\rho$ is associative but $\rho$ is clearly not a group homomorphism when $\beta \neq 1$. This is a slightly unusual definition of a projective representation. Classically, a projective representation is a group homomorphism $\rho: G \to \PL(V)$. While these definitions are equivalent, for our purposes this second definition is unsatisfactory because we care about the explicit value of the $2$-cocycle $\beta$. When mapped through to $\PL(V)$, constants disappear and as such there would be no difference between different values of $\beta$.

    As an example of the difference between these perspectives, consider the case when $G = \m Z/2\m Z = \{0, 1\}$. If we view projective representations as a map to $\PL(V)$, the only one dimensional projective representation is the trivial one. However, there are non trivial one dimensional $\beta$-projective representations. For example when $\beta$ is given by $\beta(1, 1) = -1$ and all other pairs map to $1$, two projective representations are given by $1 \mapsto i$ and $1 \mapsto -i$.

    The following is brief introduction to projective character theory with the aim of building up to a proof of Theorem \ref{Group Extension Theorem} and Corollary \ref{irr:irr:correspondance}.

    \begin{mydef}
        Two projective representations $(V, \rho)$, $(W, \psi)$ are linearly equivalent if there exists an isomorphism $f: V \to W$ such that $f\circ \rho(g) = \psi(g) \circ f$ for all $g \in G$.
    \end{mydef}

    There is another coarser notion of equivalence called projective equivalence. It is essentially the same definition but the equality occurs in $\PL(W)$ instead of $\GL(W)$. The two $\beta$-projective representations of $\m Z/2\m Z$ given a moment ago are examples of representations that are projectively but not linearly equivalent.

    \begin{mydef}
        Given two $\beta$-projective representations $(V, \rho)$, $(W, \psi)$, we can form a new $\beta$-projective representation $(V\oplus W, \rho\oplus \psi)$ called the direct sum.
    \end{mydef}

    \begin{mythm}
        Let $(V, \rho)$ be a $\beta$-projective representation and $W \subset V$ a $G$ invariant subspace. Then $(V, \rho)$ can be decomposed into the direct sum of two $\beta$-projective representations, $(W, \rho|_{W})$ and $(W', \rho|_{W'})$ where $W\oplus W' = V$ and $\rho|_{W}$ represents the action of $\rho$ on $V$ restricted to $W$.
    \end{mythm}
    \begin{proof}
        See \cite{cheng:1996}.
    \end{proof}

    \begin{mydef}
        A projective representation $(V, \rho)$ is said to be irreducible if the only $G$ invariant subspaces are $V$ and $\{0\}$.
    \end{mydef}

    \begin{mythm}
        All projective representations can be decomposed into a direct sum of irreducible representations. This decomposition is unique up to reordering and linear equivalences.
    \end{mythm}
    \begin{proof}
        See \cite{cheng:1996}.
    \end{proof}

    Recall Definition \ref{proj char}, defining the character of a projective representation,
    \[
        \chi_{\rho}(g) = \tr\big(\rho(g)\big).
    \]

    \begin{mythm}
        A $\beta$-projective representation is entirely determined by its character.
    \end{mythm}
    \begin{proof}
        See \cite{cheng:1996}.
    \end{proof}

    We are now able to prove Theorem \ref{Group Extension Theorem}, restated below for convenience.  

    \begin{mythm}
        Let $G$ be a group and $\beta$ a unitary $2$-cocycle. Then there exists a group $G_{\beta}$  called the group extension of $G$ and an injective map $f$ from $\beta$-projective representations of $G$ to linear representations of $G_{\beta}$.
    \end{mythm}
    \begin{proof}
        First to construct $G_{\beta}$. As $\beta$ is unitary it takes values in the group of $|G|$'th roots of unity in $\m C$. This group is naturally identified with $A = C_{|G|} = \langle x | x^{|G|}\rangle$ by the isomorphism $e^{\frac{2\pi i}{|G|}} \mapsto x$. Let $\tilde{\beta}$ be image of $\beta$ under the isomorphism into $A$. As a set, $G_{\beta} = G \times A$ and its group structure comes from a twisted multiplication map
        \[
            (g, x^n)\times (h, x^m) = (gh, \tilde{\beta}(g, h)x^{n + m}).
        \]
        Associativity of the multiplication follows from the associativity of the multiplication in $G$ and $A$ as well as the $2$-cocycle condition which $\beta$ satisfies. As $\beta$ is normalized, there is an identity element $(e, 1)$, and inverses are given by $(g, x^m)^{-1} = (g^{-1}, \big(\beta(g, g^{-1})\big)^{-1}x^{-m}$. This proves that $G_{\beta}$ is indeed a group.

        Note that $A$ is isomorphic to the normal subgroup of $G_{\beta}$ generated by $(e, x)$ and that $G$ is not a subgroup of $G_{\beta}$ but is a quotient group given by $G = G_{\beta}/A$.

        Next let $(V, \rho)$ be some $\beta$-projective representation of $G$. Define $f(\rho): G_{\beta} \to GL(V)$ by 
        \[
            f(\rho)\big(g, x^m\big) = e^{\frac{2m\pi i}{|G|}}\rho(g).
        \]
        The following calculation shows that $(V, f(\rho))$ is indeed a linear representation of $G_{\beta}$.
        \begin{align*}
            f(\rho)\big(g, x^m\big)f(\rho)\big(h, x^n\big) & = e^{\frac{2(m+n)\pi i}{|G|}}\rho(g)\rho(h)
            \\ & = \beta(g, h)e^{\frac{2(m+n)\pi i}{|G|}}\rho(gh)
            \\ & = f(\rho)\big(gh, \tilde{\beta}(g, h) x^{n + m}\big)
            \\ & = f(\rho)\big((g, m)\times (h, m)\big).
        \end{align*}
        The injectivity of $f$ is immediate from noting that $f(\rho)\big(g, 1\big) = \rho(g)$ and this completes the proof.
    \end{proof}

    Comparing this construction to the one given in Section \ref{cha:CompData}, observe that $(e, x)$ corresponds to $x$ and that for any generator $g_i$, $g_i$ corresponds to $(g_i, 1)$.

    We now move on to a proof of Corollary \ref{irr:irr:correspondance} again restated below.

    \begin{mycor}
        The $\beta$-projective representation $(V, \rho)$ will be irreducible if and only if $(V, f(\rho))$ is an irreducible linear representation.
    \end{mycor}
    \begin{proof}
        Let $W$ be a $G_{\beta}$ invariant subspace of $V$. Then as $f(\rho)\big(g, 1\big) = \rho(g)$, $W$ must also be a $G$ invariant subspace. Conversely, assume that $W$ is not a $G_{\beta}$ invariant subspace. Then there exists some element $w \in W$ and $(g, x^m) \in G_{\beta}$ such that $f(\rho)\big(g, x^m)\big)(w) \notin W$. Therefore $\rho(g)(w) = e^{\frac{-2m\pi}{|G|}} f(\rho)\big(g, x^m)(w) \notin W$ and so $W$ is not a $G$ invariant subspace. Therefore $G$ invariant subspaces of $V$ correspond to $G_{\beta}$ invariant subspace of $V$.
    \end{proof}

\section{The Modular Equivalence between \texorpdfstring{$\Rep(D^{\omega}G)$}{Rep(Dw G)} and \texorpdfstring{$\B Z(\V^{\omega^{-1}}) G$}{Z(Vec w-1 G)}}
\label{app:app2}

    We provide a detailed proof of Theorem \ref{Representation Equivalence Ribbon}.

    This theorem will be proven incrementally, through the following lemmas.

    \begin{mylem} \label{Representation Equivalence Categories}
        There is an equivalence on the level of categories between $\Rep(D^{\omega}G)$ and $\B Z\big(\V^{\omega^{-1}} G\big)$.
    \end{mylem}
    \begin{mylem} \label{Representation Equivalence Monoidal}
        The equivalence on the level of categories between $\Rep(D^{\omega}G)$ and $\B Z\big(\V^{\omega^{-1}} G\big)$ extends to a monoidal equivalence.
    \end{mylem}
    \begin{mylem} \label{Representation Equivalence Briaded}
        The monoidal equivalence between $\Rep(D^{\omega}G)$ and $\B Z\big(\V^{\omega^{-1}} G\big)$ can be made into a braided equivalence.
    \end{mylem}
    \begin{mylem} \label{Representation Equivalence Pivotal}
        The monoidal equivalence between $\Rep(D^{\omega}G)$ and $\B Z\big(\V^{\omega^{-1}} G\big)$ can be made into a pivotal equivalence.
    \end{mylem}

    The functor $Eq: \Rep(D^{\omega}G) \to \B Z\big(\V^{\omega^{-1}} G\big)$ is defined as follows. Let $(V, \rho)$ be a left $D^{\omega}G$ module. Recall that $V$ must split as
    \[
        V = \bigoplus_{g \in G} V_g
    \]
    where $V_g = \im(\rho(\delta_g\ov{e})) = \Span(v_{g_1}, \cdots, v_{g_{k_g}})$ and, given an element $\delta_g \ov{x}$, $\rho(\delta_g \ov{x})$ is an invertible linear map from $V_{x^{-1}gx} \to V_g$, acting as $0$ on $V_h$ for $h \neq x^{-1}gx$.

    This shows that a left $D^{\omega}G$ module is naturally a $G$-graded vector space.  However, in order to be an element of $\B Z\big(\V^{\omega^{-1}} G\big)$ we also require a half braiding, $\ov{\beta}$. Denote $W_g$ to be the one dimensional $G$-graded vector space sitting over the element $g \in G$ and let $\beta_g = \rho\left(\sum_{k \in G} \delta_k\ov{g}\right)$. Define a map $\ov{\beta}_{g}$ from $W_g \otimes V \to V \otimes W_g$ by $\ov{\beta}_{g} = (\beta_g \otimes 1) \circ F$. With $F$, the flip map which sends $v\otimes w \to w\otimes v$.

    \begin{mylem}
        This map $\ov{\beta}_{g}$ defines a half braiding on the simple objects of $\V^{\omega^{-1}} G$ and so extends to a half braiding on all of $\V^{\omega^{-1}} G$.
    \end{mylem}

    \begin{proof}
        First we show that $\beta_g^{-1}$ exists (So $\ov{\beta}_{g}$ is invertible) and $\ov{\beta}_{g}$ is a $G$-graded map. Then we will show that $\ov{\beta}_{g}$ is natural and that it satisfies the half braiding equality $\omega^{-1}_{V, g, h} \circ (\ov{\beta}_g \otimes 1) \circ (\omega^{-1}_{g, V, h})^{-1} \circ (1 \otimes \ov{\beta}_h) \circ \omega^{-1}_{g, h, V} = \ov{\beta}_{gh}$.

        Initially, observe that $\left(\sum_{k \in G} \delta_k\ov{g}\right)^{-1} = \sum_{l \in G} \theta_{glg^{-1}}(g, g^{-1})^{-1} \delta_l\ov{g^{-1}}$ and so $\beta_g^{-1}$ exists. Next $\beta_g$ restricted to $V_{g^{-1}h}$ is a map from $V_{g^{-1}h}$ to $\to V_{hg^{-1}}$. Hence as
        \[
            W_g \otimes V = \bigoplus_{h \in G} V_{g^{-1}h},
        \]
        and
        \[
            V \otimes W_g = \bigoplus_{h \in G} V_{hg^{-1}}
        \]        
        $\ov{\beta}_g$ is indeed a $G$-graded map.

        Naturality of $\ov{\beta}_g$ is immediate as it is acts by the identity map on $W_g$ and so we are left with showing the half braiding condition. Hence pick a representative element $(u, w, v) \in W_g \otimes W_h \otimes V$ with $v = \sum_{k \in G} v_k$ where each $v_k \in V_k$. The right hand side of the proposed equivalence acts upon this element by 
        \[
            u \otimes w \otimes v \xmapsto{\ov{\beta}_{gh}} \beta_{gh}(v) \otimes u \otimes w = \sum_{k \in G} \beta_{gh}(v_k) \otimes u \otimes w.
        \]
        In order to calculate the action of the left hand side, observe that if $v_k \in V_k$, then $\beta_h(v_k) \in V_{hkh^{-1}}$. Using this fact we can derive that left hand side will act by
        \begin{align*}
            u \otimes w \otimes v & = \sum_k u \otimes w \otimes v_k
            \\ & \xmapsto{\omega^{-1}_{g, h, V}} \sum_k \omega^{-1}(g, h, k) u \otimes w \otimes v_k
            \\ & \xmapsto{1 \otimes \ov{\beta}_h} \sum_{k} \omega^{-1}(g, h, k) u \otimes \beta_h(v_k) \otimes w
            \\ & \xmapsto{(\omega^{-1}_{g, V, h})^{-1}} \sum_{k} \frac{\omega^{-1}(g, h, k)}{\omega^{-1}(g, hkh^{-1}, h)} u \otimes \beta_h(v_k) \otimes w
            \\ & \xmapsto{\beta_g \otimes 1} \sum_{k} \frac{\omega^{-1}(g, h, k)}{\omega^{-1}(g, hkh^{-1}, h)} \beta_g(\beta_h(v_k)) \otimes u \otimes w
            \\ & \xmapsto{\omega^{-1}_{V, g, h}} \sum_{k} \frac{\omega^{-1}(g, h, k)\omega^{-1}(ghk(gh)^{-1}, g, h)}{\omega^{-1}(g, hkh^{-1}, h)} \beta_g(\beta_h(v_k)) \otimes u \otimes w
            \\ & = \sum_{k} \theta^{-1}_{ghk(gh)^{-1}}(g, h) \beta_g(\beta_h(v_k)) \otimes u \otimes w.
        \end{align*}
        Simply observe now that
        \begin{align*}
            \beta_g(\beta_h(v_k)) & = \rho\Big(\Big(\sum_{l, m} \nabla(\delta_l\ov{g}, \delta_m\ov{h}\Big), v_k\Big)
            \\ & = \rho\Big(\sum_{l} \theta_{l}(g, h)\delta_l\ov{gh}, v_k\Big)
            \\ & = \theta_{ghk(gh)^{-1}}\rho(\delta_{ghk(gh)^{-1}}\ov{gh}, v_k)
            \\ & = \theta_{ghk(gh)^{-1}} \beta_{gh}(v_k).
        \end{align*}
        Substituting this in shows that the left hand side and the right hand side act identically and so the proposed morphism is indeed a half braiding.
    \end{proof}

    This defines how $Eq$ maps objects in $D^{\omega}G$ to objects in $\B Z\big(\V^{\omega^{-1}} G\big)$ sending $(V, \rho) \mapsto (V, \ov{\beta})$. On, morphisms, $Eq$ will correspond to the identity map as a morphism that commutes with the action of $D^{\omega}G$ must commute with $\rho(\delta_g\ov{e})$ and $\beta_g$ which means it must respect the derived $G$-grading and will commute with the half braiding.

    It also needs to be shown that $Eq$ has an inverse\footnote{In general an equivalence of categories requires a slightly weaker condition than in inverse existing but in this case the constructed functor $Eq$ is really on the nose invertible.}. For brevity I will only give a brief outline of this inverse as explaining the finer details is essentially a rehash of the work done above.

    Given, an object $(V, \ov{\beta})$ in $\B Z\big(\V^{\omega^{-1}} G\big)$, the corresponding $D^{\omega}G$ module is created by forgetting the $G$ grading on $V$ and equipping it with the action
    \[
        \rho(\delta_g\ov{x}) = (\ov{\beta}_x) \circ P(x^{-1}gx)
    \]
    where $P(h)$ is the projection of $V$ onto $V_h$. Morphisms will again be mapped across by the identity map.

    This proves Lemma \ref{Representation Equivalence Categories}.

    Now we show that this equivalence between $\Rep\big(D^{\omega}G\big)$ and $\B Z\big(\V^{\omega^{-1}} G\big)$ is a monoidal one.

    \begin{proof} [Proof of Lemma \ref{Representation Equivalence Monoidal}]
        Let $(V, \rho)$ and $(W, \psi)$ be two left $D^{\omega}G$ modules. Denote the associated elements of $\B Z\big(\V^{\omega^{-1}} G\big)$ by $Eq(V, \rho) = (V, \ov{\alpha})$ and $Eq(W, \psi) = (W, \ov{\beta})$. Initially we show that
        \[
            Eq(V, \rho) \otimes Eq(W, \psi) = Eq((V, \rho) \otimes (W, \psi)) = Eq(V \otimes W, (\rho \otimes \psi)\circ \Delta).
        \]
        Recall that
        \[
            \Big(V \otimes W\Big)_g = \bigoplus_{hk = g} V_h \otimes W_k = \bigoplus_{h \in G} V_h \otimes W_{h^{-1}g}
        \]
        and $V_g$ is exactly the fixed subspace of $V$ corresponding to the projection by the action of $\delta_g\ov{e}$. Then the action on $V\otimes W$ of $\delta_g\ov{e}$ is given by
        \[
            (\rho \otimes \psi)\big(\Delta(\delta_g\ov{e})\big) = (\rho \otimes \psi)\Big(\sum_{h \in G}\delta_h\ov{e}\otimes\delta_{h^{-1}g}  \ov{e}\Big) = \Big(\sum_{h \in G}\rho(\delta_h\ov{e})\otimes\psi(\delta_{h^{-1}g}\ov{e})\Big).
        \]
        Hence as required, $(\rho \otimes \psi)\big(\Delta(\delta_g\ov{e})\big)$ is exactly a projection onto the space
        \[
            \bigoplus_{h \in G} V_h \otimes W_{h^{-1}g}
        \]
        and so the monoidal structure on the vector space part matches up.

        Next recall that the tensor product on half braidings is given on a simple object $U_g \in \V^{\omega^{-1}} G$ by
        \[
            (\ov{\alpha} \otimes \ov{\beta})_g = (\omega^{-1}_{V, W, g})^{-1} \circ (1 \otimes \ov{\beta}_g) \circ \omega^{-1}_{V, g, W} \circ (\ov{\alpha}_g \otimes 1) \circ (\omega^{-1}_{g, V, W})^{-1}.
        \]
        Similarly, the half braiding from the tensor product of $D^{\omega}G$ modules is
        \[
            \ov{\mu}_g = (\mu_g \otimes 1) \circ F
        \]
        where
        \begin{align*}
            \mu_g & = (\rho \otimes \psi) \circ \Delta\left(\sum_{h \in G} \delta_h\ov{g}\right)
            \\ & = (\rho \otimes \psi)\left(\sum_{h, k \in G} \gamma_g(k, k^{-1}h)\delta_k\ov{g}\otimes \delta_{k^{-1}h}\ov{g}\right)
            \\ & = \sum_{h, k \in G} \gamma_g(k, k^{-1}h)\rho(\delta_k\ov{g}) \otimes \psi(\delta_{k^{-1}h}\ov{g})
            \\ & = \sum_{l, m \in G} \gamma_g(l, m)\rho(\delta_l\ov{g}) \otimes \psi(\delta_m\ov{g}).
        \end{align*}
        Hence, consider these acting on a sample element $(u \otimes v \otimes w) \in U_g \otimes V\otimes W$. Additionally, assume that $v$ and $w$ split as $v = \sum_{g} v_g$ and $w = \sum_{g} w_g$ with $v_h \in V_h$ and $w_k \in W_k$. Then
        \begin{align*}
            u \otimes v \otimes w & = \sum_{h, k} u\otimes v_h \otimes w_k
            \\ & \xmapsto{(\omega^{-1}_{g, V, W})^{-1}} \sum_{h, k} \omega(g, h, k) u \otimes v_h \otimes w_k
            \\ & \xmapsto{\ov{\alpha}_g \otimes 1} \sum_{h, k} \omega(g, h, k) \alpha_g(v_h) \otimes u \otimes w_k
            \\ & \xmapsto{\omega^{-1}_{V, g, W}} \sum_{h, k} \frac{\omega(g, h, k)}{\omega_{ghg^{-1}, g, k}} \alpha_g(v_h) \otimes u \otimes w_k
            \\ & \xmapsto{1 \otimes \ov{\beta}_g} \sum_{h, k} \frac{\omega(g, h, k)}{\omega_{ghg^{-1}, g, k}} \alpha_g(v_h) \otimes \beta_g(w_k) \otimes u
            \\ & \xmapsto{(\omega^{-1}_{V, W, g})^{-1}} \sum_{h, k} \frac{\omega_{g, h, k}\omega_{ghg^{-1}, gkg^{-1}, g}}{\omega_{ghg^{-1}, g, k}} \alpha_g(v_h) \otimes \beta_g(w_k) \otimes u
            \\ & = \sum_{h, k, l, m} \gamma_{g}(ghg^{-1}, gkg^{-1}) \rho(\delta_l\ov{g}, v_h) \otimes \psi(\delta_{m}\ov{g}, w_k) \otimes u
            \\ & = \sum_{h, k, l, m} \gamma_{g}(ghg^{-1}, gkg^{-1}) \rho(\delta_l\ov{g}, v_h) \otimes \psi(\delta_{m}\ov{g}, w_k) \otimes u
            \\ & = \sum_{l, m} \gamma_{g}(l, m) \rho(\delta_l\ov{g}, v_{g^{-1}lg}) \otimes \psi(\delta_{m}\ov{g}, w_{g^{-1}kh}) \otimes u
            \\ & = \mu_g(v \otimes w) \otimes u
            \\ & = \ov{\mu}_g(u \otimes v \otimes w).
        \end{align*}
        Therefore the equivalence preserves the operation of tensor product on objects. It is easy to show that the operation of tensor product on morphisms is also preserved. It remains to show that $Eq$ preserves the rest of the monoidal structure.

        Note that the identity object of $\Rep(D^{\omega}G)$ is given by $(\m{C}, \epsilon)$ and it can be easily checked that $F(\m{C}, \epsilon) = (V_e, 1)$ where $1$ represents the trivial braiding, which is the identity object of $\B Z\big(\V^{\omega^{-1}} G\big)$. In addition, observe that both the unitors of both categories are trivial and the associator\footnote{This comes from the quasi-associative strcucture on $D^{\omega} G$ and must be inverted} on $\Rep(D^{\omega}G)$, $\Phi^{-1}$ exactly corresponds to multiplying an element in $U_g\otimes V_h \otimes W_k$ by $\omega^{-1}(g, h, k)$. This is exactly the associator on $\B Z\big(\V^{\omega^{-1}} G\big)$ and so with a trivial tensorator and unit isomorphism the previously constructed categorical equivalence becomes a monoidal equivalence. 
    \end{proof}

    Next we show that this equivalence extends to the braiding. 

    \begin{proof} [Proof of Lemma \ref{Representation Equivalence Briaded}]
        Let $(V, \rho)$ and $(W, \psi)$ be two left $D^{\omega}G$ modules with corresponding objects in $\B Z\big(\V^{\omega^{-1}} G\big)$, $(V, \ov{\alpha})$ and $(W, \ov{\beta})$.

        Recall that the braiding we have given to $\B Z\big(\V^{\omega^{-1}} G\big)$ is induced from the half braiding to be
        \[
            \sigma_{V, W} = \ov{\beta}_V.
        \]
        On $\Rep(D^{\omega}G)$, the braiding comes from the quasitriangular element $R$, and is given by $F \circ (\rho \otimes \psi)(R)$. Hence we compute
        \begin{align*}
            \sigma_{V, W} = \ov{\beta}_V & = (\beta_V \otimes 1) \circ F 
            \\ & = \Bigg(\sum_{g \in G} \psi\Big(\sum_{h \in G} \delta_h\ov{g}\Big) \otimes 1\Bigg) \circ F
            \\ & = \Bigg(\sum_{g, h \in G} \psi\big(\delta_h\ov{g}\big) \otimes \rho\big(\delta_h\ov{e}\big)\Bigg) \circ F
            \\ & = F \circ \Bigg(\sum_{g, h \in G} \rho\big(\delta_h\ov{e}\big) \otimes \psi\big(\delta_h\ov{g}\big) \Bigg)
            \\ & = F \circ (\rho \otimes \psi) (R).
        \end{align*}
        Therefore, $\Rep(D^{\omega}G)$ is braided equivalent to $\B Z\big(\V^{\omega^{-1}} G\big)$.
    \end{proof}

    Finally we need to show that the rigid and pivotal structures are also equivalent under this map.
    \begin{proof} [Proof of Lemma \ref{Representation Equivalence Pivotal}]
        There are two parts to this theorem. Initially, we must show that the duals with their evaluation and co-evaluation maps are equivalent. Then we can move on the the pivotal structure.

        Dual objects must be preserved by the equivalence as being duals is a property of a pair of object. We just need to show then that our choices for evaluation and co-evalutaion maps also align. Recall that the two constants $\alpha$ and $\beta$ from the definition of a quasi Hopf algebra are given as $\alpha = \textbf{1}$ and $\beta = \sum_{g}\omega(g, g^{-1}, g)\delta_g\ov{e}$. The means that the evaluation map is exactly function application and the co-evaluation map is given by $(1 \otimes \rho^*(\beta)) \circ \eta_V$ where, picking any basis $v_i$ for $V$, $\eta_V$ is the map
        \[
            1 \xmapsto{\eta_V} \sum_{i} v_i\otimes v^i.
        \]
        Observe that $S(\beta) = \beta^{-1} = \sum_{g}\omega^{-1}(g, g^{-1}, g)\delta_g\ov{e}$. Hence, this $\rho^*(\beta)$ term will exactly multiply the $g$-th graded piece by $\omega^{-1}(g, g^{-1}, g)$ and so this exactly corresponds to the rigid structure for right duals on $\B Z(\V^{\omega^{-1}} G)$. Therefore the equivalence preserves the right rigid structure. We will deal with the left rigid structure after showing that the equivalence preserves the pivotal isomorphism.

        To find the pivotal structure on $\Rep(D^{\omega}G)$ we first make the observation that for all $h \in D^{\omega}G$,
        \[
            S^2(h) = \beta^{-1} h \beta.
        \]
        This means that $\beta^{-1}h = S^2(h)\beta^{-1}$. Letting $\ov{\phi}$ be the usual pivotal structure on $\V$ and choosing an arbitrary $f \in V^*$ we find that
        \begin{align*}
            \ov{\phi} \circ \rho(\beta^{-1}) \circ \rho(h)(v)(f) & = \rho(\beta^{-1}h, v)^{**}(f)
            \\ & = f\Big(\rho\big(S^2(h)\beta^{-1}, v\big)\Big)
            \\ & = \rho^*(S(h), f)\big(\rho(\beta^{-1}, v)\big)
            \\ & = \big(\rho(\beta^{-1}, v)\big)^**\Big(\rho^*(S(h), f)\Big)
            \\ & = \rho^{**}\Big(h, \big(\rho(\beta^{-1}, v)\big)^{**}\Big)(f)
            \\ & = \rho^{**}(h) \circ \ov{\phi} \circ \rho(\beta^{-1})(v)(f).
        \end{align*}
        Hence, the pivotal structure on $\Rep(D^{\omega}G)$ is given by $\ov{\phi} \circ \rho(\beta^{-1})$. As $\beta^{-1} = \sum_{g}\omega(g, g^{-1}, g)^{-1}\delta_g\ov{e}$ this is exactly multiplying the component of $v$ in $V_g$ by $\omega^{-1}(g, g^{-1}, g)$ and which is precisely the the pivotal structure on $\B Z\big(\V^{\omega^{-1}}G\big)$ and so clearly the equivalence extends to to these pivotal structures.

        Recall that the pivotal isomorphism gives a canonical isomorphism between the left and right duals. Hence, as the equivalence preserves the right rigid structure and and the pivotal isomorphism it must also preserve the left rigid structure.
    \end{proof}

    Over the course of these four lemma's we have shown that the equivalence between $\B Z\big(\V^{\omega^{-1}}G\big)$ and $\Rep(D^{\omega}G)$ is a braided pivotal equivalence. As equivalences of categories will also preserve properties, Theorem \ref{Representation Equivalence Ribbon} immediately follows and thus $\B Z\big(\V^{\omega^{-1}}G\big)$ and $\Rep(D^{\omega}G)$ are equivalent ribbon fusion categories.

\section*{Acknowledgments}
We would like to thank Corey Jones and the reviewer for their helpful comments. Scott Morrison was supported by ARC grants DP160103479 and FT170100019

\bibliographystyle{unsrtnat}
\bibliography{paper}

\end{document}